\documentclass[12pt]{amsart}
\usepackage{amsfonts, amsmath, amsthm}
\usepackage{epsfig}
\usepackage{psfrag}

\newtheorem{pro}{Proposition}[subsection]
\newtheorem{thm}[pro]{Theorem}
\newtheorem{lem}[pro]{Lemma}
\newtheorem{clm}[pro]{Claim}

\newtheorem{cor}[pro]{Corollary}

\theoremstyle{definition}
\newtheorem{dfn}[pro]{Definition}
\newtheorem{ex}[pro]{Example}

\theoremstyle{remark}

\title{2-Normal Surfaces} 
\date{\today}
\address{Mathematics Department, University of Illinois at Chicago}
\email{bachman@math.uic.edu}
\author{David Bachman}
\begin{document}

\begin{abstract}
We define a 2-normal surface to be one which intersects every 3-simplex of a triangulated 3-manifold in normal triangles and quadrilaterals, with one or two exceptions. The possible exceptions are a pair of octagons, a pair of unknotted tubes, an octagon and a tube, or a 12-gon. 

In this paper we use the theory of critical surfaces developed in \cite{crit} to prove the existence of topologically interesting 2-normal surfaces. Our main results are (1) if a ball with normal boundary in a triangulated 3-manifold contains two almost normal 2-spheres then it contains a 2-normal 2-sphere and (2) in a non-Haken 3-manifold with a given triangulation the minimal genus common stabilization of any pair of strongly irreducible Heegaard splittings can be isotoped to an almost normal or a 2-normal surface.
\end{abstract}
\maketitle

\noindent
Keywords: Heegaard Splitting, Stabilization, Normal Surface, Minimal Surface.

\section{Introduction.}

In \cite{jr:88} Jaco and Rubinstein pioneered the viewpoint that in triangulated 3-manifolds the normal surfaces of Kneser \cite{kneser:29} are analogous to stable minimal surfaces. Rubinstein pushed this idea further in \cite{rubinstein:93} by introducing {\it almost normal surfaces}, the analogue of an unstable minimal surface of index one. Guided by this analogy he showed that any strongly irreducible Heegaard splitting can be isotoped to an almost normal surface. This was a crucial step in showing that in a non-Haken 3-manifold there are at most a finite number of distinct Heegaard splittings of a given genus \cite{jr:03}, resolving a conjecture of Waldhausen. 


In this paper we present a theory of {\it 2-normal} surfaces. A 2-normal surface is defined to be one which intersects every 3-simplex of a triangulated 3-manifold in normal triangles and quadrilaterals, with one or two exceptions. The possible exceptions are a pair of octagons, a pair of unknotted tubes, an octagon and a tube, or a 12-gon. (See Theorem \ref{t:classify} for a more precise definition.) This theory is guided by the viewpoint that 2-normal surfaces are analogous to unstable minimal surfaces of index two. 

We develop the theory of 2-normal surfaces by examining the interplay between ``critical" Heegaard splittings and triangulations. {\it Criticality} was introduced in \cite{crit} as a useful combinatorial condition on the compressing disks of a Heegaard splitting, much like the condition of {\it strong irreducibility} \cite{cg:87}. In \cite{crit} we show the following:

\begin{thm} 
\label{t:nocrit}
Suppose $M$ is an irreducible 3-manifold with no closed incompressible surfaces and at most one Heegaard splitting (up to isotopy) of each genus. Then $M$ does not contain a critical Heegaard splitting. 
\end{thm}

This shows that critical Heegaard splittings are a non-trivial class of surfaces. What shows that they are an interesting class is the following:

\begin{thm}
\label{t:mgcs}
Suppose $M$ is a small 3-manifold whose boundary, if non-empty, is incompressible and $F$ and $F'$ are distinct strongly irreducible Heegaard splittings of $M$ which induce the same partition of $\partial M$. Then the minimal genus common stabilization of $F$ and $F'$ is critical.
\end{thm}

We begin this paper by defining what it means for a surface to be critical {\it relative to a 1-manifold} (Section \ref{s:Relative}) and exploring some of the useful properties of such surfaces. When the 1-manifold in question is a suitable subset of the 1-skeleton of a triangulation we show that such surfaces are 2-normal (Section \ref{s:PL}). In the remainder of the paper (Section \ref{s:exist}) we show that a critical Heegaard surface can either be made critical or strongly irreducible relative to a given 1-manifold. By combining this machinery with Theorem \ref{t:mgcs} we obtain the main result of this paper:

\medskip
\noindent {\bf Theorem \ref{t:main2}} {\it The minimal genus common stabilization of any pair of strongly irreducible Heegaard splittings of a non-Haken 3-manifold is isotopic to an almost normal or 2-normal surface in any triangulation.}
\medskip

At present it is not known whether or not there is an algorithm which will take two Heegaard splittings and decide how many times one must stabilize one to obtain a stabilization of the other. Any normal surface-type algorithm to do this would have to enumerate all possibilities for the minimal genus common stabilization of the two. The significance of Theorem \ref{t:main2} is that it provides a framework for such an enumeration.

In developing the machinery necessary to prove Theorem \ref{t:main2} we also prove the following interesting result:

\medskip
\noindent {\bf Theorem \ref{t:almostnormalpair}} {\it Let $M$ be a 3-manifold with traingulation $T$ and $B$ an embedded 3-ball in $M$ with normal boundary. Suppose $B$ contains no vertices of the 0-skeleton of $T$ and the only normal 2-spheres in $B$ are copies of $\partial B$. If $B$ contains two almost normal 2-spheres, each with an octagon, then $B$ contains a 2-normal 2-sphere with either two octagons or one 12-gon.}
\medskip


This theorem asserts the existence of surfaces similar to those conjectured by Rubinstein in his program to show that any action of $\mathbb Z^n$ on $S^3$ must be standard \cite{rubinstein:96}. Note that the existence of an almost normal 2-sphere with octagon in any submanifold $B$ that satisfies the hypotheses is an important step in Rubinstein's algorithm to recognize the 3-sphere \cite{rubinstein:93} (see also \cite{thompson:94}).

Theorems \ref{t:almostnormalpair} and \ref{t:main2} are further evidence that almost normal and 2-normal surfaces are the appropriate analogues of unstable minimal surfaces of indices one and two. One would expect that if one had distinct index one minimal surfaces then there should be an index two minimal surface. The analogue of this is precisely the assertion of the above theorems.

\tableofcontents

\subsection{Basic Definitions.}
In this section we give definitions of some of the standard terms that will be used throughout the paper. The expert in 3-manifold theory can easily skip this. 

A 2-sphere in a 3-manifold which does not bound a 3-ball on either side is called {\it essential}. If a manifold does not contain an essential 2-sphere, then it is referred to as {\it irreducible}.

A loop on a surface is called {\it essential} if it does not bound a disk in the surface. Given a surface $F$ in a 3-manifold $M$ a {\it compressing disk} for $F$ is a disk $D \subset M$ such that $F \cap D=\partial D$ and such that $\partial D$ is essential on $F$. If we let $D \times I$ denote a thickening of $D$ in $M$ then to {\it compress $F$ along $D$} is to remove $(\partial D) \times I$ from $F$, and replace it with $D \times \partial I$. A surface for which there are no compressing disks is called {\it incompressible}. A manifold which contains an incompressible surface or an essential 2-sphere is {\it Haken}. 

A 3-manifold is {\it small} if either it is closed and non-Haken, or it is has incompressible boundary, and every incompressible surface is boundary parallel. 

A {\it compression body} is a 3-manifold which can be obtained by starting with some surface, $F$, forming the product, $F \times I$, attaching some number of 2-handles to $F \times \{1\}$, and capping off any remaining 2-sphere boundary components with 3-balls. The boundary component, $F \times \{0\}$, is often referred to as $\partial _+$. A {\it Heegaard splitting} of a 3-manifold $M$ is an expression of $M$ as a union $W \cup _F W'$, where $W$ and $W'$ are compression bodies that intersect in $F=\partial _+W=\partial _+ W'$. Such a splitting is {\it nontrivial} if neither $W$ nor $W'$ are products. If $W \cup _F W'$ is a Heegaard splitting of $M$ then we say $F$ is a {\it Heegaard surface}. A Heegaard surface $F$ is {\it strongly irreducible} if every compressing disk for $F$ in $W$ intersects every compressing disk for $F$ in $W'$. A {\it stabilization} of $F$ is a new Heegaard surface which is the connect sum of the standard genus 1 Heegaard surface in $S^3$ and $F$. Another way to define a stabilization is by ``tunneling" a 1-handle out of $W$ and attaching it to $\partial _+ W'$. If this is done in such a way so as to make the definition symmetric in $W$ and $W'$ then one arrives at a stabilization. The Riedemeister-Singer theorem states that given any two Heegaard surfaces, $F$ and $F'$, there is always a stabilization of $F$ which is isotopic to a stabilization of $F'$.

\section{Critical Surfaces}

\subsection{Definitions.}

In this section we summarize the basic definitions of \cite{crit}. 

Let $F$ be an embedded surface separating a 3-manifold $M$ into a ``red" and a ``blue" side. If $C$ and $C'$ are compressing disks for $F$ then we say $C$ is equivalent to $C'$ if there is an isotopy of $M$ taking $F$ to $F$, and $C$ to $C'$ (we do allow $C$ and $C'$ to be on opposite sides of $F$). We denote the equivalence class of a compressing disk $C$ as $[C]$. 

We now define a 1-complex $\Gamma (F)$. For each equivalence class of compressing disk for $F$ there is a vertex of $\Gamma (F)$. Two (not necessarily distinct) vertices are connected by an edge if there are representatives of the corresponding equivalence classes on opposite sides of $F$ which intersect in at most a point. 

\begin{dfn}
A vertex of $\Gamma (F)$ is {\it isolated} if it is not the endpoint of any edge.
\end{dfn}

\begin{dfn}
If we remove the isolated vertices from $\Gamma (F)$ and are left with a disconnected 1-complex then we say $F$ is {\it critical}.
\end{dfn}

Equivalently, $F$ is critical if there exist two edges of $\Gamma(F)$ that can not be connected by a 1-chain. 

Suppose $D$ and $E$ are compressing disks on opposite sides of $F$ such that $|D \cap E| \le 1$. Then we denote the edge of $\Gamma (F)$ which connects $[D]$ to $[E]$ as $D-E$. If $[D]=[D']$ then we write $D \sim D'$. Hence, a chain of edges in $\Gamma (F)$ may look something like
\[D_1-E_2 \sim E_3-D_4 \sim D_5 - E_6 -D_7\]
Many of the proofs in \cite{crit} and the present paper follow by producing such chains. 

If $C$ and $C'$ are compressing disks for $F$ such that $[C]$ and $[C']$ are in the same component of $\Gamma(F)$ then we say $C$ is {\it $\Gamma$-path connected} to $C'$.

\subsection{Local properties of critical surfaces.}
\label{s:local}

In this section we prove results about critical surfaces that will be needed in the remainder of the paper. However, as these results deal only with critical surfaces this section should really be regarded as an extension of \cite{crit}.

\begin{lem}
\label{l:3sphere}
$S^3$ does not contain any critical surfaces.
\end{lem}

\begin{proof}
Let $F$ be a critical surface in $S^3$. Compress $F$ as much as possible to the red side. If $F$ compresses to a collection of spheres then the red side is a handlebody. If not, then it compresses down to a surface $F_r$. Similarly, compress $F$ as much as possible to the blue side. If the blue side is not a handlebody then $F$ compresses to a surface $F_b$. If both the red and blue sides are handlebodies then $F$ is a critical Heegaard splitting of $S^3$ and we can apply Theorem \ref{t:nocrit}. 

Let $M$ denote the submanifold of $S^3$ cobounded by $F_b$ and $F_r$ (If one of these surfaces is empty then let $M$ denote the side of the other that contains $F$). By construction, every compressing disk for $\partial M$ must lie in $M$. Also by construction, $F$ is a critical Heegaard splitting for $M$ so Corollary 5.11 of \cite{crit} implies that $\partial M$ is incompressible in $M$. We conclude then that $\partial M$ is incompressible in $S^3$, a contradiction. 
\end{proof}

\begin{cor}
\label{c:3ball}
Let $F$ be a connected critical surface in a 3-manifold $M$ and $B$ be an embedded 3-ball. Then $B$ cannot contain $F$. 
\end{cor}

\begin{proof}
Suppose $F \subset B$. Any curve that bounds a compressing disk for $F$ in $M$ also bounds a compressing disk for $F$ in $B$. Hence, not only is $F$ critical in $M$ but it is also critical in $B$. Now, glue to $B$ another ball, forcing $F$ to be a critical surface in $S^3$. This contradicts Lemma \ref{l:3sphere}.
\end{proof}

\section{Relatively Critical Surfaces}
\label{s:Relative}

\subsection{Definitions.}
\label{s:relK}

We now repeat many of our definitions with respect to some properly embedded 1-manifold $K$. Let $M^K$ denote $M$ with a regular neighborhood of $K$ removed. If $X$ is any subset of $M$ then let $X^K=X \cap M^K$. Let $F \subset M$ be an embedded, closed, separating surface. Let $D$ be an embedded disk in $M$ such that $\partial D=\alpha \cup \beta$, where $D \cap F=\alpha$ and $D \cap K=\beta$. Then $D$ will be referred to as a {\it relative compressing disk} for $F$. In general, a {\it $K$-compression} for $F$ is any relative compressing disk for $F$ or any compressing disk for $F^K$. 

Suppose $F$ separates $M$ into a ``red" and a ``blue" side. The term {\it red disk} will be used to refer to any $K$-compression on the red side of $F$ and {\it blue disk} for any $K$-compression on the blue side. 

{\bf Notation.} Red disks will usually be denoted with the letter ``$D$" (usually with some subscript) and blue disks with the letter ``$E$".

If $C$ and $C'$ are $K$-compressions for $F$ then we say $C$ is equivalent to $C'$ if there is an isotopy of $M$ taking $F$ to $F$ and $C$ to $C'$, which fixes $K$ setwise (we do allow $C$ and $C'$ to be on opposite sides of $F$). 

We now define the 1-complex $\Gamma (F;K)$ in a similar manner as before. For each equivalence class of $K$-compression for $F$ there is a vertex of $\Gamma (F;K)$. Suppose $D$ is a red disk and $E$ is a blue disk. Then $D-E$ is an edge of $\Gamma (F;K)$ if:
\begin{enumerate}
    \item $D$ and $E$ are disjoint, or
    \item $D$ and $E$ are compressing disks such that $|D \cap E| =1$, or 
    \item $D$ and $E$ are relative compressing disks such that $|D \cap E| =1$ and $D \cap E \in K \cap F$.
\end{enumerate}

\begin{dfn}
$F$ is {\it strongly irreducible relative to $K$} if there are $K$-compressions on opposite sides of $F$ but $\Gamma (F;K)$ contains no edges.
\end{dfn}

\begin{ex}
In \cite{bachman:98} we show that if a knot is in both {\it thin position} (see Examples \ref{e:thin} and \ref{e:thin2} below) and {\it bridge position} (i.e. there is a thick level which separates the minima of $K$ from the maxima) then a bridge sphere is strongly irreducible relative to $K$. 
\end{ex}

We now come to our main definition:

\begin{dfn}
If we remove the isolated vertices from $\Gamma (F;K)$ and are left with a disconnected 1-complex then we say $F$ is a {\it critical surface relative to $K$}. When it is understood what $K$ is we simply call $F$ {\it relatively critical}. 
\end{dfn}

\subsection{Local properties of relatively critical surfaces}

\begin{cor}
\label{c:rel3ball}
Let $F$ be a connected surface which is critical relative to some properly embedded 1-manifold $K$. Let $B$ be a 3-ball embedded in the complement of $K$. Then $B$ cannot contain $F$. 
\end{cor}

\begin{proof}
Note that if $F$ is a critical surface relative to $K$, but $F \cap K=\emptyset$, then $F$ is also just a critical surface in the complement of $K$. Hence, by Corollary \ref{c:3ball} no relatively critical surface can lie entirely inside a 3-ball which is embedded in the complement of $K$.
\end{proof}

\begin{lem}
\label{l:rellocal}
Let $K$ be a properly embedded 1-manifold in a 3-manifold $M$ such that $M^K$ is irreducible. Let $F$ be a critical surface relative to $K$. Suppose $D_1$ and $E_2$ are a red and blue disk such that $[D_1]$ and $[E_2]$ lie in different components of $\Gamma (F)$, but neither are isolated. Then there does not exist a ball $B$ embedded in the complement of $K$ such that $D_1$ and $E_2$ are contained in the interior of $B$.
\end{lem}

\begin{proof}
Let $\Lambda$ denote the set of embedded balls in $M$ such that for each $B \in \Lambda$
\begin{enumerate}
    \item there is a disk $D_1' \subset int(B)$ which is $\Gamma$-path connected to $D_1$,
    \item there is a disk $E_2' \subset int(B)$ which is $\Gamma$-path connected to $E_2$, and
    \item every loop of $F \cap \partial B$ is essential on $F$.
\end{enumerate}
If the Lemma is false then clearly there is a ball for which conditions 1 and 2 hold. By assumption, the complement of $K$ is irreducible. Hence, an easy innermost disk argument shows that such a ball may be isotoped to satisfy condition 3 as well. We proceed then under the assumption that $\Lambda$ is non-empty.

Let $B$ be an element of $\Lambda$. Note that the closure of every disk component of $\partial B-F$ is either a red or a blue disk for $F$. Let $\Lambda _b$ denote the subset of $\Lambda$ such that for each $B \in \Lambda _b$ there is a blue disk among the closure of the components of $\partial B-F$. Similarly, let $\Lambda _r$ denote the subset of $\Lambda$ such that for each $B \in \Lambda _r$ there is a red disk among the closure of the components of $\partial B-F$

\begin{clm}
\label{c:union}
$\Lambda=\Lambda _b \cup \Lambda _r$.
\end{clm}

\begin{proof}
If the claim is false then there is a ball $B \in \Lambda$ such that $F \cap \partial B=\emptyset$. As the interior of $B$ contains a red disk by assumption, $F$ must lie entirely inside $B$. This contradicts Corollary \ref{c:rel3ball}.
\end{proof}

\begin{clm}
\label{c:disjoint}
$\Lambda_b \cap \Lambda _r=\emptyset$.
\end{clm}

\begin{proof}
Let $B$ be an element of both $\Lambda _b$ and $\Lambda _r$. Then among the closures of the components of $\partial B-F$ there is a blue disk $E$ and a red disk $D$. By definition of $\Lambda$ there are also disks $D_1'$ and $E_2'$ inside $B$ which are $\Gamma$-path connected to $D_1$ and $E_2$. As $D$ and $E$ are subsets of $\partial B$ they must be disjoint from $D_1'$ and $E_2'$. If $\partial D \ne \partial E$ then $D$ and $E$ are disjoint. Otherwise, we may push one of them slightly out of $B$ to make them disjoint. In either case we end up with the following contradictory chain: 
\[D'_1-E-D-E'_2\]
\end{proof}

Our assumption that $\Lambda \ne \emptyset$, along with Claim \ref{c:union}, implies that at least one of $\Lambda _b$ and $\Lambda _r$ is non-empty. Henceforth we will assume that $\Lambda_b \ne \emptyset$. Our proof will be symmetric otherwise. 

\begin{clm}
\label{c:nonempty}
Let $D_2$ be any red disk which is $\Gamma$-path connected to $E_2$. Then $D_2 \cap \partial B \ne \emptyset$ for any $B \in \Lambda _b$. 
\end{clm}

\begin{proof}
As $B \in \Lambda$ there is a red disk $D_1' \subset int(B)$ which is $\Gamma$-path connected to $D_1$. As $B \in \Lambda _b$ there is a blue disk $E$ among the closures of the components of $\partial B -F$. Hence, $D_1' \cap E=\emptyset$. If $D_2 \cap E = \emptyset$ as well then we would have the chain $D_1'-E-D_2$, a contradiction. As $E \subset \partial B$ it must be that $D_2 \cap \partial B \ne \emptyset$. 
\end{proof}

We are now prepared to state our minimality assumption. Consider all triples of the form $(B,E_2', D_2)$, where $B \in \Lambda _b$, $E_2'$ is a blue disk in the interior of $B$ which is $\Gamma$-path connected to $E_2$, and $D_2$ is a red disk such that $E_2'-D_2$ is an edge of $\Gamma(F)$. Our assumption that $[E_2]$ was not an isolated vertex guarantees the existence of $D_2$. Henceforth, we will assume that $(B,E_2', D_2)$ has been chosen among all such triples so that $|D_2 \cap \partial B|$ is minimal. 

Note first that $D_2 \cap \partial B$ does not contain any loops, as all such curves can be removed by a standard innermost disk argument. Let $A$ denote the closure of a component of $D_2-\partial B$ such that $A \cap \partial B$ is a single arc $\alpha$. Such an arc exists by Claim \ref{c:nonempty}. Note that since there are always at least two choices for $A$, and $|D_2 \cap E_2'| \le 1$, we may assume that $A \cap E_2'=\emptyset$.

\begin{clm}
There is a red disk $D_1' \subset int(B)$ which is $\Gamma$-path connected to $D_1$ such that $D_1' \cap A=\emptyset$. 
\end{clm}

\begin{proof}
As $B \in \Lambda$ there is a disk $D_1'' \subset int(B)$ which is $\Gamma$-path connected to $D_1$. Furthermore, since $B \in \Lambda _b$ there is a blue disk $E$ among the closures of the components of $\partial B -F$.

We will recursively define a sequence of compressing disks $\{D^i\}$ for $F$ which terminates with a disk with the desired properties. Like $D_1''$, each disk we will construct will be contained in the interior of $B$ and will be $\Gamma$-path connected to $D_1$. If at any stage we construct a disk $D^i$ such that $D^i \cap A = \emptyset$ then stop and let $D_1'=D^i$. Begin by letting $D^1=D_1''$. 

Let $\beta$ be an arc of $D^{i-1} \cap A$ which is outermost on $A$. Then $\beta$ cuts off the subdisk $A'$ of $A$ and divides $D^{i-1}$ into two disks, $C$ and $C'$. At least one of $A' \cup C$ and $A' \cup C'$ is a compressing disk for $F$, since $D^{i-1}$ was. Choose one that is and denote this disk as $D^{i}$. Now, since both $D^{i-1}$ and $A'$ were contained in the interior of $B$ we have $D^{i} \subset int(B)$. Hence, $D^{i} \cap E=\emptyset$ and we have the following chain in $\Gamma(F)$:
\[D^{i}-E-D_1''\]
Furthermore, note that $|D^{i} \cap A|<|D^{i-1} \cap A|$, so our construction must terminate with a disk which is disjoint from $A$. 
\end{proof}

We now use the disk $A$ to guide an isotopy of $\partial B$, removing one arc of $D_2 \cap \partial B$, and transforming $B$ to the ball $B'$. As $A \cap D_1'=\emptyset$ and $A \cap E_2'=\emptyset$ we still have $D_1' \subset int(B')$ and $E_2' \subset int(B)$. Hence, the only way in which we have not contradicted the minimality of $|D_2 \cap \partial B|$ is if $B' \notin \Lambda _b$. 

For any $B \in \Lambda$ there are at least two disk components of $\partial B-F$. Claim \ref{c:disjoint} implies that for our particular choice of $B$ the closure of all such disks are blue. Let $D$ and $D'$ then represent two such blue disks. Recall that $A \cap \partial B=\alpha$. As $B' \notin \Lambda _b$ it must be that $\alpha \cap D$ and $\alpha \cap D'$ are non-empty. Hence, $\alpha$ must be an arc on $\partial B$ which connects $\partial D$ to $\partial D'$. The effect of the isotopy guided by $A$ is to create a tangency between $D$ and $D'$ along $\alpha$ which subsequently resolves to create a single disk. As both $D$ and $D'$ were blue this new disk must also be blue. Hence, by definition $B'$ must be an element of $\Lambda _b$, a contradiction.
\end{proof}

\subsection{Spanning Surfaces}
\label{s:spanning}

Our goal in this section is to examine the interplay between relative critical surfaces and spanning surfaces of knots and links. What we are after is the following:

\begin{lem}
\label{l:index2Gabai}
Let $K$ be a properly embedded 1-manifold in a 3-manifold $M$, let $F$ be a relatively critical surface, and let $S$ be a spanning surface for $K$. Then we can isotope $F$ so that $S$ does not contain any red or blue disks for $F$ and so that no loop of intersection of $S \cap F$ is inessential on both surfaces.
\end{lem}

{\bf Note.} Some readers may recognize this as the ``index 2" analogue of a corresponding ``index 1" theorem of Gabai (Lemma 4.4 of \cite{gabai:87}). Gabai's result is that if $K$ is a knot in thin position in $S^3$ (see Examples \ref{e:thin} and \ref{e:thin2} for relevant definitions) and $S$ is a spanning surface for $K$ then each thick level $F$ can be isotoped so that $S$ does not contain any red or blue disks for $F$ and so that no loop of intersection of $S \cap F$ is inessential on both surfaces.

\begin{proof}
The proof follows that of Theorem 5.1 of \cite{crit} almost verbatim, with the caveat that the terms {\it red disk} and {\it blue disk} have slightly more general meanings now. As this proof is a significant portion of that paper we will only sketch the argument here and describe the necessary modifications. 

The original proof is in four stages. In the first, a map $\Phi: S \times D^2 \rightarrow M$ is constructed. In the second, we break up $D^2$ into regions and use $\Phi$ to label each. The labelling is defined so that the existence of an unlabelled region would imply the Lemma. In Stage 3 we define a 2-complex $\Pi$ such that $H_1(\Pi)$ is non-trivial and use our labelling of $D^2$ to define a map from $D^2$ to $\Pi$ (assuming there are no unlabelled regions). Finally, in Stage 4 we show that our map from $D^2$ to $\Pi$, when restricted to $\partial D^2$, is non-trivial on homology. As $D^2$ is contractable this provides a contradiction. 

The only thing that needs to change for us is Stage 1 of the above argument. The rest of the argument follows verbatim. In Stage 1 we begin by letting $D_0-E_0$ and $D_1-E_1$ be edges in different components of $\Gamma (F)$ such that $D_0 \cap E_0=D_1 \cap E_1=\emptyset$. We then prove that this can always be arranged. That is, we show that if some component of $\Gamma (F)$ contains an edge representing disks that intersect in a point ({\it i.e.} a destabilization) then it contains an edge that represents disjoint disks ({\it i.e.} a weak reduction). Unfortunately, this is not necessarily true in the relative setting. Hence, for our purposes we will let $D_0-E_0$ and $D_1-E_1$ be edges in different components of $\Gamma (F;K)$ such that if $D_i$ or $E_i$ is a compressing disk then $D_i \cap E_i=\emptyset$. If for some $i=0$ or 1, $D_i$ and $E_i$ are both relative compressions then we make no restriction on $D_i \cap E_i$, other than the usual $|D_i \cap E_i| \le 1$. 

The next step in Stage 1 is to construct a sequence of disks $\{D_{\frac {i}{n}}\}_{i=1}^{n-1}$ such that $D_{\frac {i}{n}} \cap D_{\frac {i+1}{n}} =\emptyset$ for all $i$ between 0 and $n-1$. This construction follows that of the original proof verbatim. Similarly, we construct a sequence of disks $\{E_{\frac {i}{m}}\}_{i=1}^{m-1}$ such that $E_{\frac {i}{m}} \cap E_{\frac {i+1}{m}} =\emptyset$ for all $i$ between 0 and $m-1$.

Next, we let $\{U_i\}_{i=0}^{n}$ and $\{V_i\}_{i=0}^{m}$ denote neighborhoods of the disks $\{D_{\frac {i}{n}}\}$ and $\{E_{\frac {i}{m}}\}$ which satisfy certain constraints. The {\it original} constraints are as follows:

\begin{enumerate}
    \item $U_i \cap U_{i+1}=\emptyset$, for $0 \le i <n$
    \item $V_i \cap V_{i+1}=\emptyset$, for $0 \le i <m$
    \item $U_0 \cap V_0=\emptyset$
    \item $U_n \cap V_m=\emptyset$
\end{enumerate}

Unfortunately, if $D_0$ and $E_0$ are relative compressing disks which meet in a point then condition 3 cannot be satisfied. Hence, we will insist that condition 3 be satisfied only when $D_0 \cap E_0=\emptyset$ and condition 4 when $D_1 \cap E_1 =\emptyset$. 

Finally, we define a set of isotopies. For each $i$ between 0 and $n$ we let $\gamma ^i :M \times I \rightarrow M$ be an isotopy such that $\gamma ^i _0 (x)=x$ for all $x \in M$, $\gamma ^i _t (x)=x$ for all $t$ and all $x$ outside of $U_i$, and $\gamma ^i_1(S) \cap D_{\frac {i}{n}} =\emptyset$. In other words, $\gamma ^i$ is an isotopy which pushes $S$ off of $D_{\frac {i}{n}}$ inside $U_i$. Similarly, for each $i$ between 0 and $m$ we let $\delta ^i$ be an isotopy which pushes $S$ off of $E_{\frac {i}{m}}$ inside $V_i$.

If $D_0$ (say) is a relative compressing disk then we need to describe the isotopy $\gamma^0$ in more detail. The crucial fact we need to demonstrate, in order for the remainder of the proof to go through, is that if $D_0$ and $E_0$ are relative compressing disks that meet in a point then $\gamma ^0$ and $\delta^0$ can be defined so that they commute. 

Assuming $D_0$ is a relative compressing disk that meets $E_0$ in a point we will construct the isotopy $\gamma^0$ in two phases. In the first phase we describe a way to push $S$ off of $D_0$ in a neighborhood of $K \cap D_0$. In the second we will be concerned with making $S$ disjoint from the remainder of $D_0$. 

        \begin{figure}[htbp]
        \psfrag{p}{$p$}
        \psfrag{b}{$S \cap \partial N(K)$}
        \psfrag{C}{$\partial N(K)$}
        \psfrag{S}{$S \cap D_0$}
        \vspace{0 in}
        \begin{center}
        \epsfxsize=4 in
        \epsfbox{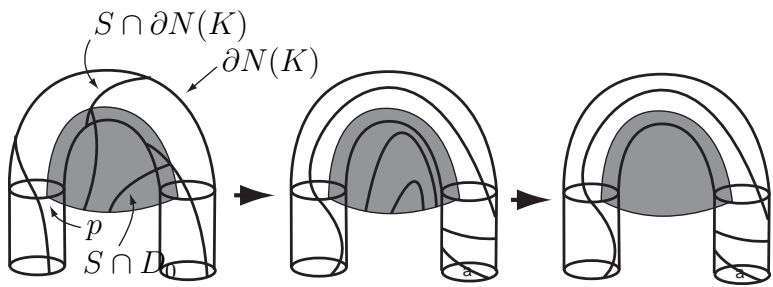}
        \caption{Defining $\gamma^0$.}
        \label{f:gamma0}
        \end{center}
        \end{figure}

To begin the first phase let $p=D_0 \cap E_0$. Now, beginning with $p$ move along $K \cap D_0$ and ``unwind" $S$, so that in the end $S \cap \partial N(K)$ is parallel to $D_0 \cap N(K)$ (see Figure \ref{f:gamma0}). This gets rid of all points of $S \cap D_0 \cap K$. In the second phase we sweep down from $K \cap D_0$ and push down any intersections of $S$ with $D_0$ that we encounter. This is also depicted in Figure \ref{f:gamma0}.

Note that this isotopy is constant near the point, $p$. Since $U_0$ and $V_0$ can be chosen so that they only overlap in a neighborhood of $p$ we may assume that $\gamma ^0$ and $\delta ^0$ commute. Similarly, if $D_1$ and $E_1$ are both relative compressions that meet in a point then we may construct  $\gamma ^n$ and $\delta ^m$ so that they commute. These are the only facts that are necessary to follow the remainder of the proof of Theorem 5.1 of \cite{crit}.
\end{proof}

\section{Triangulations}
\label{s:PL}

\subsection{Normal and Almost Normal Surfaces}

\begin{figure}[htbp]
        \vspace{0 in}
        \begin{center}
        \epsfxsize=3.25 in
        \epsfbox{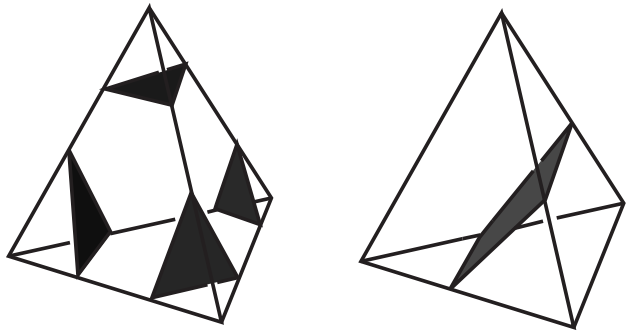}
        \caption{Normal Disks.}
        \label{f:Normal}
        \end{center}
        \end{figure}

In this section we discuss the necessary background material on normal and almost normal surfaces. A {\it normal curve} on the boundary of a tetrahedron is a simple loop which is transverse to the 1-skeleton, made up of arcs which connect distinct edges of the 1-skeleton. The {\it length} of such a curve is the number of times it crosses the 1-skeleton. A {\it normal disk} in a tetrahedron is any embedded disk whose boundary is a normal curve of length three or four, as in figure \ref{f:Normal}.

A {\it normal surface} in $M$ is the image of an embedding $p$ of some surface $F$ into $M$ such that $p(F)$ is a union of normal disks. We say $p(F)$ is an {\it almost normal surface} if it consists of all normal disks plus one additional piece in one tetrahedron. This piece can be either a disk with normal boundary of length 8 or two normal disks connected by a single unknotted tube. 

Normal surfaces were first defined by Kneser in \cite{kneser:29} and later used extensively by Haken \cite{haken:61}. Almost normal surfaces were first explored by Rubinstein in \cite{rubinstein:93} and later used by Thompson \cite{thompson:94}, Stocking \cite{stocking:96}, and the author \cite{bachman:98}, \cite{machine}.

\subsection{Intersection with the 3-skeleton}
\label{s:triangulations}

We now undertake the study of surfaces which are critical relative to a suitable subset of the 1-skeleton of a triangulation of a closed 3-manifold. Let $T$ denote such a triangulation and $T^i$ the $i$-skeleton of $T$. In this section $M$ will denote a submanifold of a 3-manifold that is bounded by a normal surface such that $T^0 \cap M=\emptyset$. $F$ will denote a critical surface of $M$ relative to $T^1 \cap M$. Recall our convention that for 1-manifolds such as $T^1 \cap M$, $M^{T^1}$ denotes the closure of the complement of a neighborhood of $T^1$ in $M$. For any subset $X$ of $M$, $X^{T^1}=X \cap M^{T^1}$.

\begin{lem}
\label{l:nlarcs}
We may isotope $F$ so that $F \cap T^2$ is a collection of normal arcs.
\end{lem}

\begin{proof}
This is a straight application of Lemma \ref{l:index2Gabai}. Here, we are considering $T^2$ as a spanning surface for $T^1$. Assume $F$ is given to us as described by Lemma \ref{l:index2Gabai}. Since $T^2$ is simply connected there can be no loops of $F \cap T^2$ which are essential on $T^2$, so Lemma \ref{l:index2Gabai} says that there are no loops at all in the intersection. If we see a non-normal arc of intersection of $F \cap T^2$ then there must be an outermost such one. Any such arc cobounds a relative compressing disk, contradicting Lemma \ref{l:index2Gabai}.   
\end{proof} 

\begin{lem}
\label{l:disjointdisks}
If $D-E$ is an edge of $\Gamma (F;T^1)$ then there is an edge $D'-E'$, in the same component as $D-E$, such that $(D' \cup E') \cap (T^2 \backslash T^1)=\emptyset$. 
\end{lem}

\begin{proof}
We begin by claiming that $D$ can be isotoped so that there is no arc $\gamma$ of $D \cap T^2$ with an endpoint on $T^1$. Note that if this is the case then $D$ must be a relative compressing disk, so that $\partial D=\alpha \cup \beta$, where $D \cap F=\alpha$, $D \cap T^1=\beta$, and $\gamma \cap T^1 \subset \beta$. Let $C$ denote the cylinder on the boundary of a neighborhood of $\beta$ which is also a subset of $\partial M^{T^1}$. Then $\beta '=D \cap C$ is an arc which  spirals around the cylinder $C$. We can now ``unwind" $\beta '$, trading points of $\beta ' \cap T^2$ for points of $\alpha \cap T^2$, as in Figure \ref{f:unwind}. 

        \begin{figure}[htbp]
        \psfrag{a}{$\alpha$}
        \psfrag{b}{$\beta '$}
        \psfrag{C}{$C$}
        \psfrag{g}{$\gamma$}
        \psfrag{T}{$T^2$}
        \vspace{0 in}
        \begin{center}
        \epsfxsize=3 in
        \epsfbox{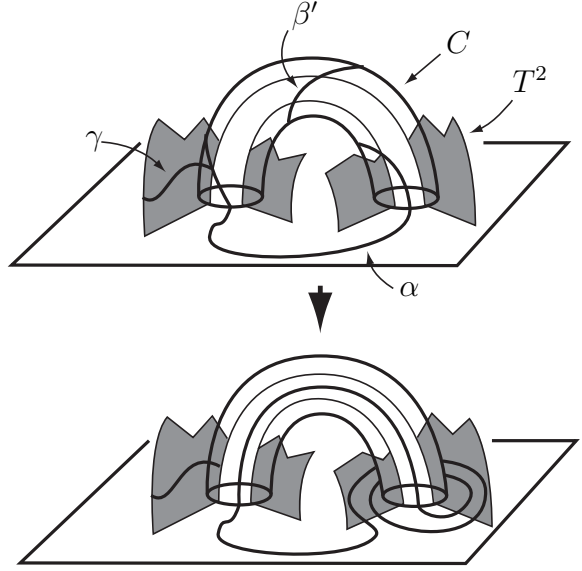}
        \caption{Unwinding $\beta'$.}
        \label{f:unwind}
        \end{center}
        \end{figure}

By similar means we can isotope $E$ so that there are no arcs of $E \cap T^2$ with endpoints on $T^1$. Suppose now the Lemma is not true. Then if $D'- E'$ is any edge in the same component of $\Gamma (F;T^1)$ as $D-E$ then $(D' \cup E') \cap (T^2 \backslash T^1) \ne \emptyset$. Choose such an edge so that
\begin{enumerate}
    \item there are no arcs of $D' \cap T^2$ or $E' \cap T^2$ with endpoints on $T^1$, and
    \item among all such disks $|(D' \cup E') \cap (T^2 \backslash T^1)|$ is minimal.
\end{enumerate}
Let $\Delta$ be a 2-simplex such that $(D' \cup E') \cap \Delta \ne \emptyset$. Note that $(\dot D' \cup \dot E') \cap \partial \Delta = \emptyset$, since the interior of any red or blue disk lies in $M^{T^1}$. There are now two cases to consider.

First, assume $D' \cap \Delta$ or $E' \cap \Delta$ contains a simple closed curve. Choose such a curve of intersection which is innermost on $\Delta$ and assume this curve lies on $D'$. So the curve cuts off a subdisk that we can use to surger $D'$. This leads to a disk $D^*$ where $D^*-E'$ is in the same component of $\Gamma (F;T^1)$ as $D-E$. But $D^* \cap (T^2 \backslash T^1)$ contains one fewer curve, contradicting the minimality of $|(D' \cup E') \cap (T^2 \backslash T^1)|$. 

The second possibility is that $(D' \cup E') \cap \Delta$ consists only of arcs contained in the interior of $\Delta$ (since we have chosen $D'$ and $E'$ so that $\partial D'$ and $\partial E'$ are disjoint from $\partial \Delta$). Suppose some such arc $\gamma$ has both endpoints on the same normal arc $n$ of $F \cap \Delta$.  Then $\gamma \cup n$ bounds a disk $U \subset \Delta$. Let $\gamma'$ be an arc of $(D' \cup E') \cap U$ which is outermost on $U$. Assume $\gamma ' \subset D'$. So $\gamma' \cup n$ bounds a subdisk $U'$ of $U$. The disk $U'$ can be used to boundary compress $D'$ in the complement of $F$. As before, this leads us to a disk $D^*$ where $|D^* \cap \Delta|<|D' \cap \Delta|$, such that $D^*-E'$ is in the same component as $D-E$, a contradiction. 

        \begin{figure}[htbp]
        \psfrag{n}{Not}
        \psfrag{o}{Outermost}
        \vspace{0 in}
        \begin{center}
        \epsfxsize=2 in
        \epsfbox{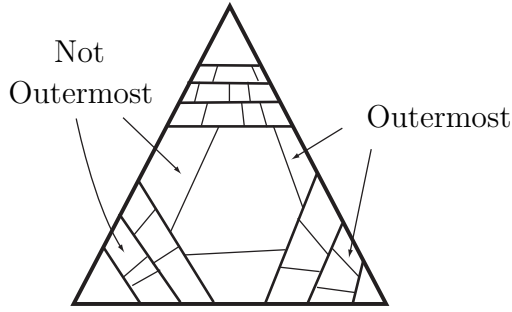}
        \caption{Defining ``outermost arcs" of $\Delta \cap (D' \cup E')$.}
        \label{f:delta}
        \end{center}
        \end{figure}

We conclude that every arc of $(D' \cup E') \cap \Delta$ must connect distinct normal arcs of $F \cap \Delta$. Choose an arc $\gamma$ of $(D' \cup E') \cap \Delta$ which is ``outermost" in the following sense: Let $n_1$ and $n_2$ be the normal arcs of $F \cap \Delta$ which $\gamma$ connects. Then there is at least one edge $e$ of $\partial \Delta$ that both $n_1$ and $n_2$ meet. The rectangle $\gamma \cup n_1 \cup n_2 \cup e$ bounds a subdisk $V$ of $\Delta$. Now, $\gamma$ is ``outermost" among the arcs of $(D' \cup E') \cap \Delta$ in the sense that $V \cap (D' \cup E')=\gamma$ (see figure \ref{f:delta}).

        \begin{figure}[htbp]
        \psfrag{e}{$e$}
        \psfrag{n}{$n_1$}
        \psfrag{o}{$n_2$}
        \psfrag{D}{$V'$}
        \psfrag{d}{$V$}
        \psfrag{g}{$\gamma$}
        \vspace{0 in}
        \begin{center}
        \epsfxsize=2 in
        \epsfbox{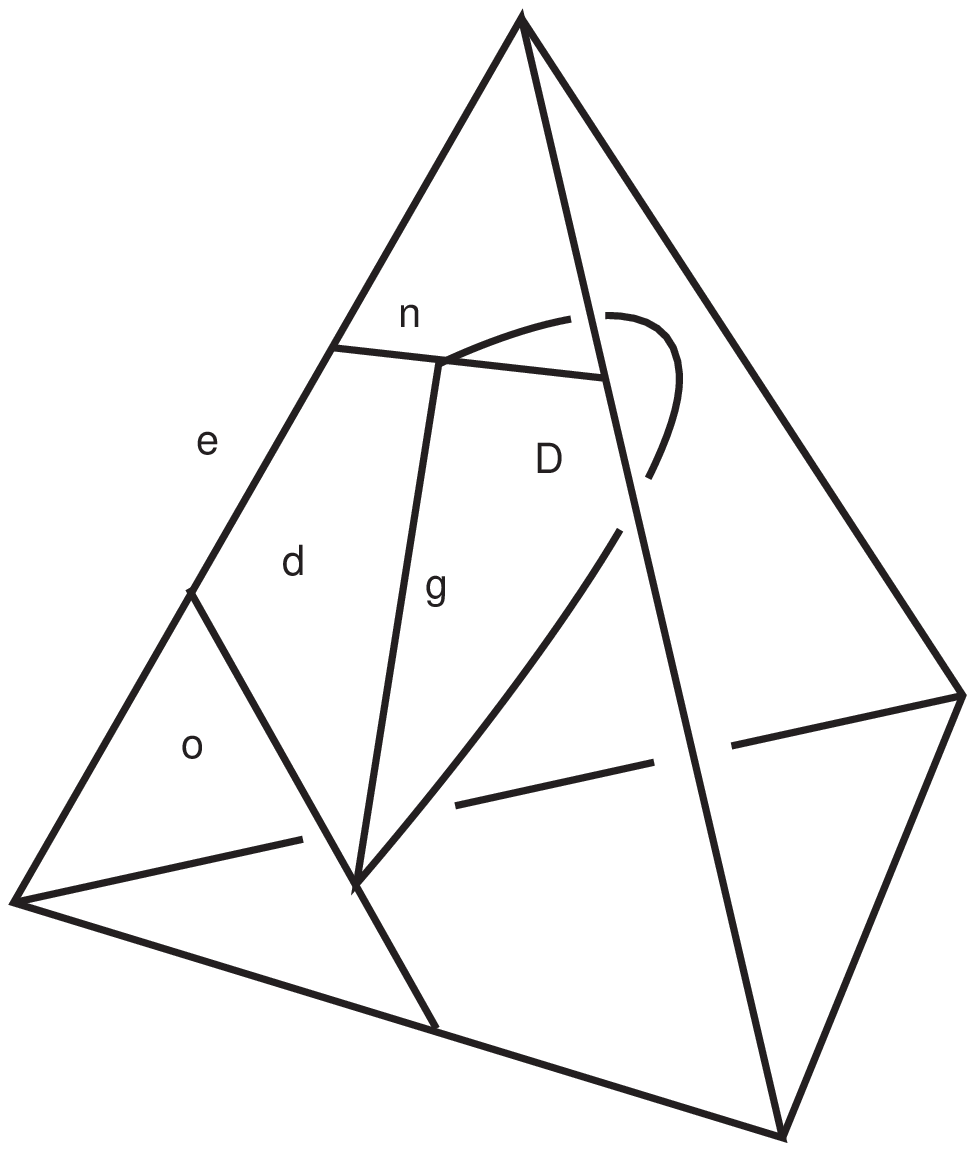}
        \caption{$D \cup D'$.}
        \label{f:edge}
        \end{center}
        \end{figure}

Suppose $\gamma \subset D'$. Note that $\gamma$ also cuts off a subdisk $V'$ of $D'$ such that $V' \cap E'=\emptyset$ and $V' \cap T^1 = \emptyset$. Let $D''$ be the disk $V' \cup V$ (see Figure \ref{f:edge}) pushed slightly off of $\Delta$. Note that since $V \cap E' =\emptyset$ (by choice of $\gamma$) and $V' \subset D'$ then $D''-E'$ is in the same component as $D-E$. But once again $|D'' \cap (T^2 \backslash T^1)|<|D' \cap (T^2 \backslash T^1)|$ which contradicts our assumption of minimality.
\end{proof}

\begin{cor}
\label{l:disjoint}
There are edges $D_1-E_1$ and $D_2-E_2$ in different components of $\Gamma (F;T^1)$ such that $D_1, D_2, E_1$, and $E_2$ are all disjoint from $T^2 \backslash T^1$.
\end{cor}

\begin{proof}
Since $F$ is relatively critical there are multiple components of $\Gamma (F;T^1)$ which contain edges. The proof is complete by two applications of Lemma \ref{l:disjointdisks}.
\end{proof}

We now produce a sequence of Lemmas which precisely narrow down all of the possible pictures of $F$ inside a tetrahedron. For the remainder of this section we assume $D_1, D_2, E_1$, and $E_2$ are as in Corollary \ref{l:disjoint}.

\begin{lem}
\label{l:topology}
The closure of every component of $F \backslash T^2$ is a disk, with the possible exception of exactly one of the following:
    \begin{enumerate}
    \item A pair of disks, connected by an unknotted tube. 
    \item Two pairs of disks, each connected by an unknotted tube, where these tubes are not nested. 
    \item Three disks connected by two unknotted tubes. This forms a pair of pants, one leg of which may be ``inside-out." 
    \item A disk with an unknotted tube which connects the disk to itself, forming a punctured torus. 
    \end{enumerate}
with the disks $D_1$, $D_2$, $E_1$, and $E_2$ depicted as in Figures \ref{f:top1} and \ref{f:top2}.
\end{lem}

{\bf Note:} Inspection of Figures \ref{f:top1} and \ref{f:top2} shows that for every compressing disk for a component of $F \backslash T^2$ there is a ``dual" relative compressing disk on the opposite side which meets the compressing disk exactly once. We will make use of this fact later. 

\begin{proof}
If, for every tetrahedron $\tau$, $F \cap \tau$ is a collection of disks then the Lemma is true. Hence, we may assume that there is a compressing disk for some component of $F \backslash T^2$.

Suppose first that there is a blue compressing disk $C$ for $F \backslash T^2$ which is disjoint from $D_i$, for $i=1$ or $2$. Then $C-D_i$ is an edge of $\Gamma (F;T^1)$. Henceforth, we rename $C$ as $E_i$. Similarly, if there is a red compressing disk for $F \backslash T^2$ which is disjoint from $E_i$ then we may choose $D_i$ to be such a disk. In other words, if we can then we choose any of the disks $D_i$ or $E_i$ to be compressing disks for $F \backslash T^2$ then we do so. So for example, if $D_1$ is a relative compressing disk then we may assume that $E_1$ meets every red compressing disk for $F \backslash T^2$. 

\medskip

\noindent \underline{Case 1.} {\it $D_1$ and $E_1$ are relative compressing disks.}

By the preceding remarks we may assume that every blue compressing disk for $F \backslash T^2$ intersects $D_1$ and every red one intersects $E_1$. Since $D_1-E_1$ is an edge of $\Gamma (F;T^1)$, $D_1$ and $E_1$ must be disjoint or meet at a point of $T^1$. Now, let $\tau$ be the tetrahedron which contains $D_1$ and $\tau '$ the tetrahedron which contains $E_1$ (it may be that $\tau =\tau '$). Let $B$ and $B'$ denote $\tau$ and $\tau '$ with a small enough neighborhood of $T^2$ removed so that the topology of $F \cap B$ and $F \cap B'$ is the same as that of $F \cap \tau$ and $F \cap \tau'$. We do this only to guarantee that $\partial B$ and $\partial B'$ are embedded in $M^{T^1}$. 

By assumption, we now have that any compressing disk for $F \backslash T^2$ must lie in either $B$ or $B'$. Now, we can use $D_1$ to guide an isotopy of $F$ (inside $B$ this just looks like a boundary compression of $F$) and then use $E_1$ to guide an isotopy. Afterward, there can be no compressing disks left inside $B$ or $B'$ since we are assuming that any such disk met either $D_1$ or $E_1$. Hence, we are left with an incompressible surface inside $B$ and $B'$. Such a surface must be a collection of disks. 

To reconstruct $F$ inside $B$ and $B'$ we begin with a collection of disks and undo a boundary compression on the blue side and then undo a boundary compression on the red side. A priori the possibilities are: 
\begin{enumerate}
    \item  a disk (Figure \ref{f:top1} (a)) 
    \item two (possibly nested) unknotted tubes ({\it i.e.} annuli) (Figure \ref{f:top1} (b)) 
    \item a single unknotted tube (Figure \ref{f:top1} (c)) 
    \item three disks connected by two unknotted tubes, forming a pair of pants with one ``leg" turned ``inside-out" (Figure \ref{f:top1} (d)) 
    \item and a disk with an unknotted tube which connects the disk to itself (Figure \ref{f:top1} (e)). 
\end{enumerate}
    
    However, in the case of nested tubes we do not see the relative compressions $D_1$ and $E_1$ at the same time, so cannot have such a configuration.  

        \begin{figure}[htbp]
        \psfrag{a}{(a)}
        \psfrag{b}{(b)}
        \psfrag{c}{(c)}
        \psfrag{d}{(d)}
        \psfrag{e}{(e)}
        \psfrag{D}{$D_1$}
        \psfrag{E}{$E_1$}
        \vspace{0 in}
        \begin{center}
        \epsfxsize=3 in
        \epsfbox{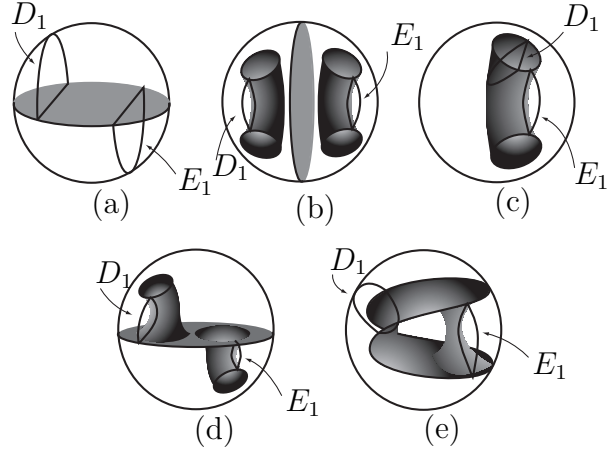}
        \caption{$D_1$ and $E_1$ are relative compressions.}
        \label{f:top1}
        \end{center}
        \end{figure}

\medskip

\noindent \underline{Case 2.} {\it $D_2$ and $E_2$ are relative compressing disks.} 

This case is symmetric with Case 1. 

\bigskip

Note that in Case 1 it was irrelevant what kinds of $T^1$-compressions ({\it i.e.} compressing disks or relative compressing disks) $D_2$ and $E_2$ are. As $D_2 \cap E_1 \ne \emptyset$ it must have been the case that $D_2 \cap B' \ne \emptyset$. Hence, in Case 1 we had completely described the topology of the component of $F \backslash T^2$ for which $D_2$ was a $T_1$-compression. Similarly, in Case 1 we had described the topology of the component of $F \backslash T^2$ for which $E_2$ was a $T_1$-compression. In other words, Case 1 and Case 2 can really be treated as being mutually exclusive, since the invocation of one obviates any invocation of the other. If we cannot invoke Case 1 or Case 2 then we move on to Case 3.

\medskip

\noindent \underline{Case 3.} {\it Either $D_1$ or $E_1$ is a compressing disk and either $D_2$ or $E_2$ is a compressing disk.} By Lemma \ref{l:rellocal} $D_1$ and $E_2$ cannot both be compressing disks and $D_2$ and $E_1$ cannot both be compressing disks. So we are left with the possibilities that $D_1$ and $D_2$ are compressing disks and $E_1$ and $E_2$ are relative compressing disks, or vice versa. Assume the former. 

Let $\tau$ now be the tetrahedron which contains $D_1$ and $E_2$. As in Case 1 let $B$ denote $\tau$ with a small enough neighborhood of $\partial \tau$ removed so that the topology of $F \cap B$ is the same as that of $F \cap \tau$. 

We now use $E_1$ to guide an isotopy of $F$ and then use $E_2$ to guide an isotopy, resulting in a surface $F'$. 

\medskip
\underline{Subcase 3.1.} {\it Some component of $F' \cap B$ is compressible.} Any compressing disk for $F'$ must be blue, since it was disjoint from both $E_1$ and $E_2$. Let $E$ denote such a disk and note that $E$ is also a compressing disk for $F$.

Let $\gamma$ be an innermost loop of $F \cap \partial B$, which is essential on $F$. The loop $\gamma$ bounds a subdisk $C$ of $\partial B$ which is a compressing disk for $F^{T^1}$. Note that $C$ is a subset of $\partial B$, which is very close to $T^2$. As $D_1$ and $D_2$ are compressing disks for $F \backslash T^2$, it must be the case that $C \cap D_i=\emptyset$ for $i=1$ and $2$. Hence, $C$ must be red. Since $E$ lies in the interior of $B$, $\partial C \cap \partial E=\emptyset$. Hence, $C-E$ is an edge of $\Gamma (F;T^1)$. 

By Lemma \ref{l:rellocal} $E$ must be $\Gamma$-path connected to $D_1$. If $D_2 \subset B$ then Lemma \ref{l:rellocal} would also say that $E$ was $\Gamma$-path connected to $D_2$, a contradiction. But if $D_2$ is not in $B$ then $E \cap D_2=\emptyset$, also implying that $E$ was $\Gamma$-path connected to $D_2$ (since $E-D_2$ would be an edge of $\Gamma (F;T^1)$). However, as $D_1$ and $D_2$ are in different components of $\Gamma (F;T^1)$ they cannot both be $\Gamma$-path connected to $E$.

\medskip
\underline{Subcase 3.2.} {\it $F' \cap B$ is incompressible.} Inside $B$, $F'$ is then a set of disks. Hence, we can reconstruct $F \backslash T^2$ (and hence, $F \cap \tau$) by starting with a collection of disks, and undoing at most two boundary compressions on the blue side. If $D_1$ and $D_2$ are both in $\tau$ then the possibilities are: 
\begin{enumerate}
    \item a disk (which cannot happen, since $D_1$ is a compressing disk for $F \cap B$) 
    \item a single unknotted annulus (Which cannot happen, since $D_1$ and $D_2$ are different compressing disks) 
    \item two unknotted annuli (Figure \ref{f:top2} (a))
    \item and a pair of pants (where neither ``leg" may be ``inside-out") (Figure \ref{f:top2} (b)). 
\end{enumerate}
    
It may be that $D_1$ and $D_2$ (and hence $E_1$ and $E_2$) are in different tetrahedra. In this case inside $B$ we only see one relative compression before we are left with a collection of disks. The only possibility for a non-disk component of $F \cap B$ in this case is a single unknotted tube. Similarly, in the tetrahedron which contained $D_2$ there may be an unknotted tube (Figure \ref{f:top2} (c)). 

        \begin{figure}[htbp]
        \psfrag{a}{(a)}
        \psfrag{b}{(b)}
        \psfrag{c}{(c)}
        \psfrag{E}{$E_1$}
        \psfrag{R}{$E_2$}
        \psfrag{D}{$D_2$}
        \psfrag{F}{$D_1$}
        \vspace{0 in}
        \begin{center}
        \epsfxsize=3 in
        \epsfbox{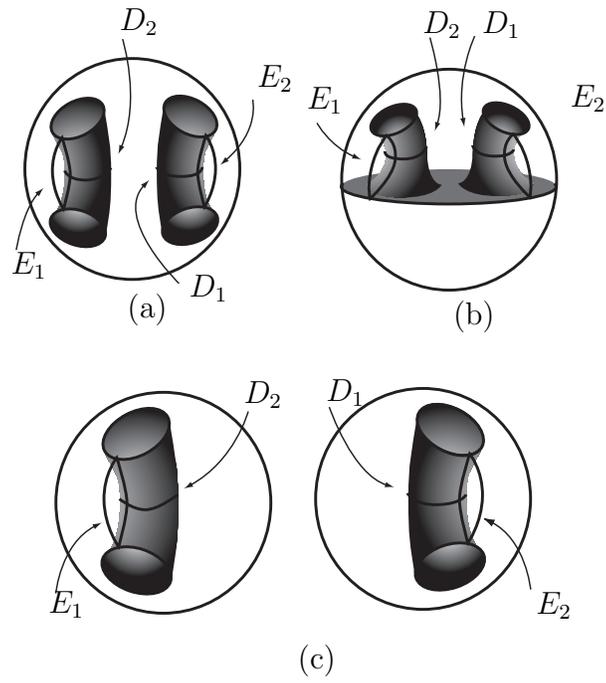}
        \caption{$E_1$ and $E_2$ are relative compressions, and $D_1$ and $D_2$ are compressions.}
        \label{f:top2}
        \end{center}
        \end{figure}

\medskip
To complete the proof of the Lemma we must now show that every component of $F \backslash T^2$ besides the exceptional ones found above is a disk. But notice that every other component is disjoint from $D_1, D_2, E_1$, and $E_2$. Hence, if some other component were compressible then there would be a red or blue disk which was disjoint from all four of these disks, a contradiction.
\end{proof}

\subsection{Intersection with the 2-skeleton}

In this section we narrow down the possibilities for the lengths of loops of $F \cap \partial \tau$, for tetrahedra $\tau$.

\begin{lem}
\label{l:tubedtriangles}
If there is a tetrahedron $\tau$ such that $F \cap \partial \tau$ contains a loop of length 8 or larger then $F \cap \tau$ does not contain an annular component whose boundary components are not normally parallel and have length 3. 
\end{lem}

        \begin{figure}[htbp]
        \vspace{0 in}
        \begin{center}
        \epsfxsize=2.5 in
        \epsfbox{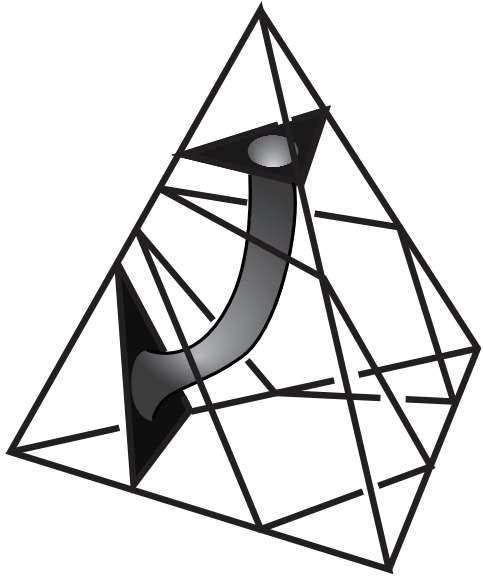}
        \caption{A tube between non-parallel normal triangles, and a curve of length 8.}
        \label{f:tubedtri}
        \end{center}
        \end{figure}

\begin{proof}
Two loops of length 3 on the boundary of a tetrahedron that are not normally parallel both intersect exactly one edge $e$ of the 1-skeleton. If these loops cobound an annular component $S$ then Lemma \ref{l:topology} implies that there is a ``dual" relative compressing disk which intersects a compressing disk for $A$ in a point and runs along $e$. But every loop of length 8 or greater meets every edge where there might be such a relative compressing disk, preventing the existence of such a disk (see Figure \ref{f:tubedtri}.) 
\end{proof}

Before proceeding any further we present two well known facts about normal loops on the boundary of a tetrahedron:

\medskip

\noindent {\bf Fact 1.} The length of a normal loop on the boundary of a tetrahedron is equal to 3 or $4n$. 

\medskip

\noindent {\bf Fact 2.} In any collection of disjoint normal loops on the boundary of a tetrahedron, if there are loops of length $4n$ then there are no loops of length $4m$ for $m \ne n$. 

\medskip

\begin{lem}
\label{l:nonparallel8}
There does not exists a tetrahedron $\tau$ such that $F \cap \partial \tau$ contains multiple loops of length 8 or greater. 
\end{lem}

\begin{proof}
If there are such loops then they are all normally parallel. Let $\alpha _1$ denote such a loop which is innermost among all such loops on $\partial \tau$, $\alpha _2$ the  closest normally parallel loop to $\alpha _1$, and $\alpha _3$ the next closest loop. Let $\beta _1$ and $\beta _2$ be normal loops of length 3 on $\partial \tau$ which are on the opposite side of $\alpha _1$ as $\alpha _2$. Finally, let $D(\alpha _i)$ be a disk in $\tau$ with boundary $\alpha _i$ and $D(\beta _i)$ a disk with boundary $\beta _i$. If $D_1$ and $D_2$ are disks then let $D_1 \# D_2$ denote the annulus obtained by connecting $D_1$ to $D_2$ by an unknotted tube. Lemmas \ref{l:topology} and \ref{l:tubedtriangles} now imply that these are the possibilities for the non-disk component of $F \cap \tau$ which contains $\alpha _1$ (up to swaping the labels of $\beta _1$ and $\beta _2$):
    \begin{itemize}
    \item $D(\alpha _1) \# D(\alpha _2)$ 
    \item $D(\alpha _1) \# D(\beta _1)$
    \item $D(\alpha _1) \# D(\beta _1) \# D(\beta _2)$ (see Figure \ref{f:pants})
    \item $D(\alpha _1) \# D(\beta _1) \# D(\alpha _2)$
    \item $D(\alpha _1) \# D(\alpha _2) \# D(\alpha _3)$
    \item Any disk connected to itself by a tube.  
    \end{itemize}
    
        \begin{figure}[htbp]
        \psfrag{A}{$\alpha _1$}
        \psfrag{a}{$\alpha _2$}
        \psfrag{B}{$\beta _1$}
        \psfrag{b}{$\beta _2$}
        \vspace{0 in}
        \begin{center}
        \epsfxsize=2 in
        \epsfbox{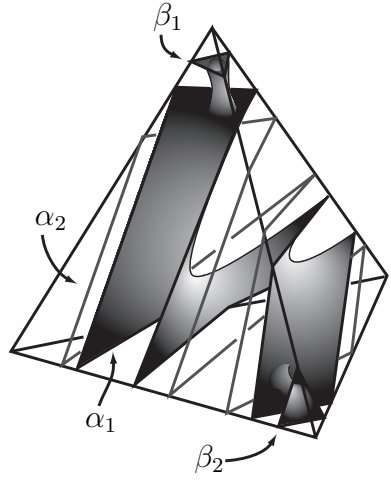}
        \caption{$D(\alpha _1) \# D(\beta _1) \# D(\beta _2)$.}
        \label{f:pants}
        \end{center}
        \end{figure}
        
In all of these cases one can find a relative compressing disk for $D(\alpha _1)$ which does not meet any other $T^1$-compression and is disjoint from $T^2$. This violates the fact that every $T^1$-compression must meet at least one of $D_1, D_2, E_1$, and $E_2$.
\end{proof}

\begin{lem}
\label{l:disknormal}
If $F \backslash T^2$ contains one of the following:
    \begin{enumerate}
    \item an $n$-gon, with $n>8$,
    \item 2 8-gons,
    \item an annular component with a boundary component of length larger than 4,
    \item a punctured torus,
    \item an 8-gon and an annular component,
    \item two annular components, or
    \item a pair of pants,
    \end{enumerate}
then every other component of $F \backslash T^2$ is a disk with boundary length 3 or 4.
\end{lem}

\begin{proof}
In all cases listed one can find a red disk $D$ for the specified piece and a blue disk $E$ such that $D-E$ is an edge of $\Gamma (F;T^1)$. (The only obstructions to the existence of these disks are ruled out in Lemmas \ref{l:tubedtriangles} and \ref{l:nonparallel8}.) As $F$ is relatively critical there is a red disk $D'$ in some other component of $\Gamma (F;T^1)$. Lemma \ref{l:disjointdisks} implies that we may choose $D'$ to be in $M \backslash T^2$. We know $E \cap D' \ne \emptyset$, since $E$ is not $\Gamma$-path connected to $D'$. 

Suppose $E^*$ is a blue disk for some other component of $F \backslash T^2$. Since $E^* \cap D =\emptyset$ we know $D - E^*$ is an edge of $\Gamma (F;T^1)$. Hence, $E^* \cap D' \ne \emptyset$, since otherwise $D-E^*-D'$ would be a chain. We now have a contradiction. The disk $D'$ cannot meet both $E$ and $E^*$ since they are $T^1$-compressions for different components of $F \backslash T^2$. 

Similarly, we can rule out any red disk for any component of $F \backslash T^2$ other than the one assumed by the Lemma. We conclude there can be no $T^1$-compressions for any other component of  $F \backslash T^2$. The absence of any compressing disks implies that each such component is a disk and the absence of relative compressing disks implies that each one has boundary of length 3 or 4.
\end{proof}

\begin{lem}
\label{l:less20}
$F \backslash T^2$ does not contain an $n$-gon such that $n > 12$.
\end{lem}

\begin{proof}
Let $U$ be a disk such that $|\partial U| > 16$, where $U$ is contained in the tetrahedron, $\tau$. (The case where $|\partial U| = 16$ will be treated separately.) By Fact 1 we know that $|\partial U|=4n$, for some integer $n$. So in fact $|\partial U|$ is at least 20. For each disk whose boundary has length greater than 4 there are arcs $\alpha _1$ and $\alpha _2$ on $\partial \tau$ with the following properties (see Figure \ref{f:20gon}): 

\begin{itemize}
    \item They connect distinct vertices of $\tau$. 
    \item The interior of each arc misses the vertices of $\tau$.
    \item For each 2-simplex $\Delta$ of $\tau$ no component of $\alpha _i \cap \Delta$ is an arc which runs from some edge back to itself. 
    \item Each arc lies in a different component of $\partial \tau \backslash \partial U$. 
    \item There are maps $\Phi _i:D^2 \times I \rightarrow \tau$ such that $\Phi _i(D^2,0)=U$, $\Phi _i(D^2,1)=\alpha _i$, and $\Phi _i | _{(D^2,[0,1))}$ is an embedding.  
\end{itemize}

\bigskip

        \begin{figure}[htbp]
        \psfrag{a}{$\alpha _1$}
        \psfrag{1}{$p_1^1$}
        \psfrag{3}{$p_1^3$}
        \psfrag{5}{$p_1^5$}
        \psfrag{7}{$p_1^7$}
        \vspace{0 in}
        \begin{center}
        \epsfxsize=3 in
        \epsfbox{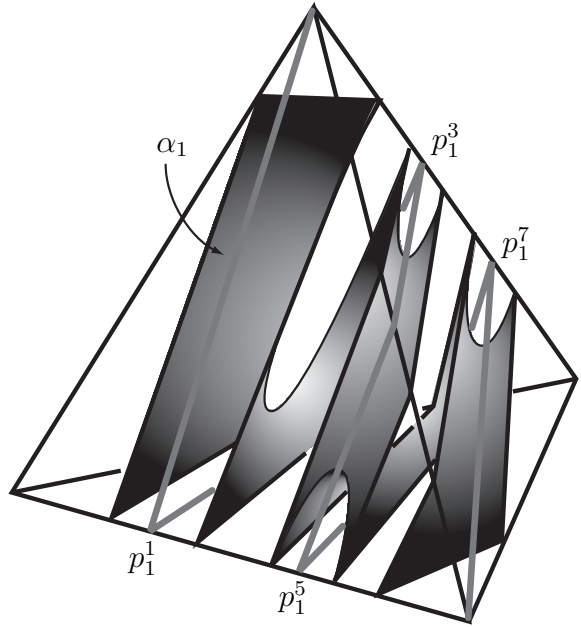}
        \caption{A 20-gon.}
        \label{f:20gon}
        \end{center}
        \end{figure}

Define a map $\phi _i$ from the image of $\Phi _i$ to $\alpha _i$ as follows: $\phi _i(\Phi _i(x,t))=\Phi _i(x,1)$. Let $p_i^j$ be the $j$th intersection point of the interior of $\alpha _i$ with $T^1$, where $1 \le j \le m$. Now, for each $(i,j)$, $\phi _i^{-1}(p_i^j)$ is a relative compressing disk for $U$. Note that the set of disks $\{\phi _1^{-1}(p_1^j)\}$ are on the opposite side of $U$ than the set $\{\phi _2^{-1}(p_2^j)\}$. We assume that the former set is on the red side. 

By Lemma \ref{l:disknormal} one of these red disks must be $D_1$ and another $D_2$. Suppose that $\phi _1 ^{-1}(p_1^k)=D_1$ and $\phi _1 ^{-1}(p_1^l)=D_2$, where $k<l$. Since $D_1 \cap E_1 =\emptyset$ and $D_2 \cap E_1 \ne \emptyset$ it follows that $\phi _1 ^{-1}(p_1^j) \cap E_1 = \emptyset$, for all $j \le k$. In particular, $\phi _1 ^{-1}(p_1^1) \cap E_1=\emptyset$ and so $\phi _1 ^{-1}(p_1^1) - E_1$ is an edge of $\Gamma (F;T^1)$. This justifies us assuming from the outset that  $\phi _1 ^{-1}(p_1^1)=D_1$. Similarly, we can show that we may assume $\phi _1 ^{-1}(p_1^m)=D_2$, $\phi _2 ^{-1}(p_2^m)=E_1$, and $\phi _2 ^{-1}(p_2^1)=E_2$.

        \begin{figure}[htbp]
        \psfrag{h}{$d_1$}
        \psfrag{k}{$d_2$}
        \psfrag{l}{$e_1$}
        \psfrag{j}{$e_2$}
        \vspace{0 in}
        \begin{center}
        \epsfxsize=2 in
        \epsfbox{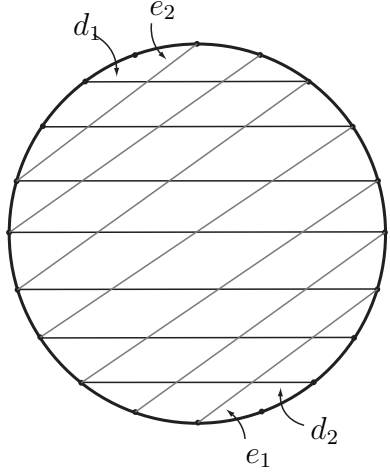}
        \caption{A 20-gon in the plane. The horizontal lines are $\{\phi _1^{-1}(p_1^j)\} \cap U$ and the slanted lines are $\{\phi _2^{-1}(p_2^j)\} \cap U$.}
        \label{f:longdisk}
        \end{center}
        \end{figure}

        \begin{figure}[htbp]
        \psfrag{h}{$d_1$}
        \psfrag{k}{$d_2$}
        \psfrag{l}{$e_1$}
        \psfrag{j}{$e_2$}
        \vspace{0 in}
        \begin{center}
        \epsfxsize=2 in
        \epsfbox{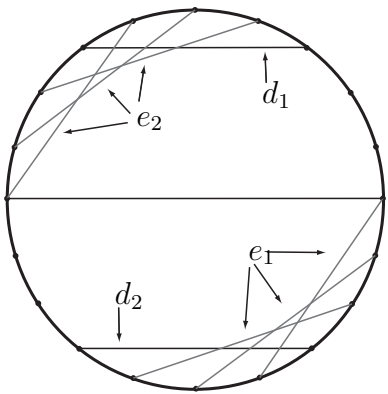}
        \caption{The three possibilities for $e_1 \cup e_2$.}
        \label{f:lows}
        \end{center}
        \end{figure}

We would now like to think of $U$ as a flat disk lying in the plane, as in Figure \ref{f:longdisk}. In this picture we have represented the arcs of $U \cap \phi _1^{-1}(\{p_1^j\})$ as horizontal lines and $U \cap \phi _2^{-1}(\{p_2^j\})$ as slanted lines. Also, there are $4n$ points marked on the boundary (with $n \ge 5$) where $\partial U \cap T^1 \ne \emptyset$. Note that all of these points are endpoints of horizontal lines except for exactly 6. Let $d_i=D_i \cap U$ and $e_i=E_i \cap U$. Now, $d_1 \cap e_2$ and $d_2 \cap e_1$ are nonempty, so there are precisely 3 possibilities for the arcs $e_1 \cup e_2$, as depicted in Figure \ref{f:lows}. Since $n \ge 5$ it follows that there are at least 5 points on each subarc of $\partial U$ between $d_1$ and $d_2$. It is now clear that there will always be a horizontal line which meets $e_1 \cup e_2$ in at most its endpoints. But this arc is part of the boundary of a red relative compressing disk which is either disjoint from $E_1$ and $E_2$ or meets them only in points of $T^1$, a contradiction.

        \begin{figure}[htbp]
        \psfrag{h}{$d_1$}
        \psfrag{k}{$d_2$}
        \psfrag{l}{$e_1$}
        \psfrag{j}{$e_2$}
        \vspace{0 in}
        \begin{center}
        \epsfxsize=2 in
        \epsfbox{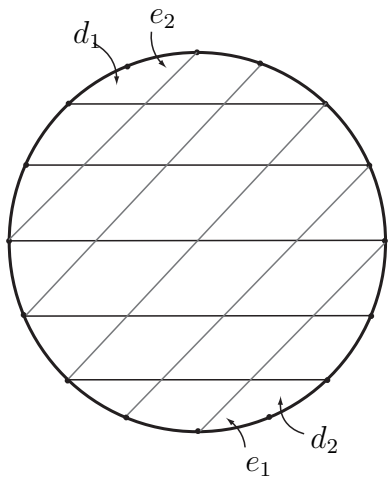}
        \caption{The unique 16-gon in the plane.}
        \label{f:flat16}
        \end{center}
        \end{figure}

Finally, to complete the proof we must rule out disks whose boundary have length 16. Such disks are special in that there is only one possibility for them in a tetrahedron, up to homeomorphisms of the tetrahedron which take edges to edges. If we repeat the argument given above we see the disk in the plane shown in Figure \ref{f:flat16}. Notice that in this picture there is a horizontal line, representing a red disk, which meets $e_1$ and $e_2$ only in its endpoints. This corresponds to a red relative compressing disk which meets both $E_1$ and $E_2$ only in a point of $T^1$, a contradiction. 
\end{proof}

\begin{lem}
\label{l:annularbound}
Every boundary component of each annular component of $F \backslash T^2$ must have length 8 or less.
\end{lem}

\begin{proof}
Suppose some annular component has a boundary component of length greater than 8. Lemma \ref{l:nonparallel8} and Fact 2 imply that the other boundary component must have length 3. Lemma \ref{l:disknormal} implies that every other component of $F \backslash T^2$ is a disk with boundary of length 3 or 4. So, $F$ is made up of normal triangles, quadrilaterals, and one exceptional piece. The exceptional piece looks like a normal triangle connected to a disk $U$ by an unknotted tube, where $|\partial U|>8$ (see Figure \ref{f:12-gon}). 

        \begin{figure}[htbp]
        \psfrag{D}{$D$}
        \psfrag{E}{$E$}
        \psfrag{e}{$E^*$}
        \vspace{0 in}
        \begin{center}
        \epsfxsize=2.5 in
        \epsfbox{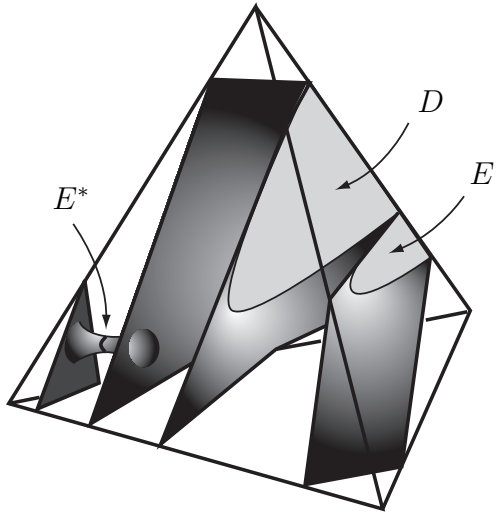}
        \caption{A 12-gon tubed to a normal triangle.}
        \label{f:12-gon}
        \end{center}
        \end{figure}

The proof is now similar to that of Lemma \ref{l:disknormal}. As $|\partial U|>8$ there are relative compressing disks $D$ and $E$ for $U$ such that $D - E$ is an edge of $\Gamma (F;T^1)$. Since $F$ is relatively critical there is some red disk $D'$ in some other component of $\Gamma (F;T^1)$. We know $E \cap D' \ne \emptyset$ since $E$ is not $\Gamma$-path connected to $D'$. 

Suppose a compressing disk $E^*$ for the annulus is blue. Note that $D \cap E^* =\emptyset$ so that $D - E^*$ is an edge of $\Gamma (F;T^1)$. Also, $E^* \cap D' \ne \emptyset$ since otherwise $D-E^*-D'$ would be a chain in $\Gamma (F;T^1)$. This is a contradiction since there is no red disk which meets both $E$ and $E^*$. 
\end{proof}

\begin{lem}
\label{l:selftube}
If there is a punctured torus component of $F \backslash T^2$ then its boundary has length 8.
\end{lem}

\begin{proof}
Suppose $P$ is a punctured torus component. Such a component has a compressing disk which, by Lemma \ref{l:topology}, has a ``dual" relative compressing disk. The boundary of this relative compressing disk must run from some edge $e$ along $P$ and back to $e$. But this implies that $P$ hits $e$ in more than one place which does not happen if $|\partial P|$ is 3 or 4. 

Lemma \ref{l:disknormal} implies that all other components are normal triangles and quadrilaterals. Suppose now that $|\partial P|>8$. Then we can think of $P$ as being obtained from a disk $U$ by attaching both ends of a tube to it. We now classify all of the red and blue disks for $P$. First, there are red and blue compressing disks $D^*$ and $E^*$ such that $|D^* \cap E^*|=1$. For each relative compressing disk of $U$ there is a relative compressing disk of $P$ which misses all compressing disks of $P$. Call such a relative compressing disk Type I. Also, for each relative compressing disk of $U$ there is a relative compressing disk of $P$ which meets either $D^*$ or $E^*$ exactly once. Call such a relative compressing disk Type II. 

Our goal now is to reach a contradiction by showing that all edges are in the same component of $\Gamma(F;T^1)$. First, since $D^*$ and $E^*$ meet once $D^*-E^*$ is an edge. Lemma \ref{l:rellocal} implies that any other edge involving a compressing disk must be in the same component as this one. Also, $D^*$ pairs with any Type I blue disk to form an edge which is obviously in the same component as $D^*-E^*$. Similarly, $E^*$ pairs with every Type I red disk to form an edge in the same component as $D^*-E^*$. We may now conclude that {\it any} edge involving a Type I disk (red or blue) is in the same component as $D^*-E^*$. 

It follows from the fact that $|\partial U|>8$ that there are relative compressions $D$ and $E$ for $U$ such that $D - E$ is an edge (see the proof of Lemma \ref{l:less20}). Let $D^I, D^{II}, E^I$ and $E^{II}$ be Type I and II red and blue disks for $F$ that correspond to $D$ and $E$. Then $D^I - E^{II}$ is an edge as well as $D^{II}- E^I$. By the remarks in the preceding paragraph these edges are in the same component as $D^*-E^*$. Now, any red Type II disk forms an edge with any blue Type II disk. In particular every red Type II disk forms an edge with $E^{II}$ and every blue Type II disk forms an edge with $D^{II}$. Hence, every Type II disk is an endpoint of an edge in the same component as $D^*-E^*$. Once again, we may conclude that {\it any} edge involving a Type II disk is in the same component as $D^*-E^*$. We conclude that all edges are in the same component, contradicting the fact that $F$ is relatively critical. 
\end{proof}

\begin{lem}
\label{l:2tubenormal}
Every boundary component of a pair of pants component of $F \backslash T^2$ has length 3 or 4. 
\end{lem}

\begin{proof}
Suppose $P$ is a pair of pants component of $F \backslash T^2$. Then $P$ can be thought of as 3 disks, $U_1, U_2$, and $U_3$, with two unknotted tubes attached. Lemma \ref{l:nonparallel8} and Fact 1 imply that at most one of these disks has boundary length greater than 4. Fact 2 implies that if one of these disks does have boundary length greater than 4 then the other two disks have boundary length 3. Suppose now that $|\partial U_1|>4$. Then $U_1$ has relative compressing disks on both sides which are disjoint from all compressing disks of $P$. (The only obstructions to this are ruled out in Lemmas \ref{l:tubedtriangles} and \ref{l:nonparallel8}.) 

Let $D$ be a compressing disk for $P$. Since $P$ is a pair of pants there is some other compressing disk which is not parallel to $D$. Either this other compressing disk or a relative compressing disk which is dual to it is blue. Let this disk be $E$. Then $D - E$ is an edge. Now let $E'$ be a blue relative compressing disk for $U_1$ which misses $D$. As in the proof of Lemma \ref{l:disknormal}, $E$ and $E'$ must meet some red disk. But inspection shows that such a red disk does not exist. We conclude that there can be no relative compressing disks for $U_1$. This implies $|\partial U_1|<8$. 
\end{proof}

If we compile the results of this section into one Theorem we have proved the following:

\begin{thm}
\label{t:classify}
Suppose $M$ is a submanifold of a 3-manifold with triangulation $T$ that is bounded by a normal surface such that $T^0 \cap M=\emptyset$. Suppose further that $F$ is a critical surface of $M$ relative to $T^1 \cap M$. Let $F$ be a critical surface of $M$ relative to $T^1$. Then in each tetrahedron $F$ looks like a collection of normal disks except for exactly one of the following:
    \begin{enumerate}
        \item Two octagons in different tetrahedra.
        \item Two unknotted tubes. These may connect two distinct pairs of normal disks or connect one normal disk to two others, forming a pair of pants.
        \item An octagon and a tube. The tube may be disjoint from the octagon, may connect the octagon to itself, or connect it to another normal disk.
        \item A 12-gon.
    \end{enumerate}
\end{thm}

We will refer to any surface described by the conclusion of Theorem \ref{t:classify} as {\it 2-normal} to indicate that such a surface fits into a classification in which there are at most 2 non-normal pieces. Following this pattern, if a surface is normal everywhere except for exactly one octagon or one unknotted tube which connects normal disks then it is tempting to refer to it as {\it 1-normal}. However, we will stick to what has become the standard terminology (coined by J.H. Rubinstein) and refer to such a surface as {\it almost normal}.

\section{The existence of relatively critical surfaces.}
\label{s:exist}

\subsection{A review of the origins of critical surfaces}
\label{s:critreview}

The goal of the rest of the paper is to answer the question ``When can a critical surface be made critical, relative to some 1-manifold?" To answer this, we must review the results from \cite{crit} which tell us precisely where critical surfaces come from. Then we must try to understand where relatively critical surfaces come from. Only then can we relate the two. 

We begin by reviewing the terminology introduced in \cite{crit}.

\begin{dfn} 
A {\it Generalized Heegaard Splitting} (GHS) of a 3-manifold $M$ is a sequence of closed, embedded, pairwise disjoint surfaces $\{F_i\} _{i=0} ^{2n}$ such that for each odd $i$ the surface $F_i$ is a non-trivial Heegaard splitting, or a union of non-trivial Heegaard splittings, of the submanifold of $M$ co-bounded by $F_{i-1}$ and $F_{i+1}$ and such that $\partial M=F_0 \amalg F_{2n}$. 
\end{dfn}

{\it Note:} We allow $F_i$ to be a union of Heegaard splittings only when the submanifold of $M$ co-bounded by $F_{i-1}$ and $F_{i+1}$ is disconnected. 

\begin{dfn}
For each odd number $i$ the surface $F_i$ of the GHS $\{F_i\}$ is said to be a {\it thick level}. Similarly, for each even value of $i$ the surface $F_i$ is said to be a {\it thin level}. 
\end{dfn}

\begin{figure}[htbp]
        \psfrag{1}{$F_0$}
        \psfrag{2}{$F_1$}
        \psfrag{3}{$F_2$}
        \psfrag{4}{$F_{2n-1}$}
        \psfrag{5}{$F_{2n}$}
        \vspace{0 in}
        \begin{center}
        \epsfxsize=1 in
        \epsfbox{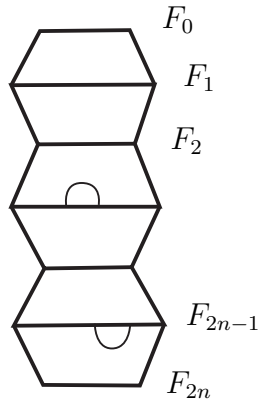}
        \caption{Schematic depicting a Generalized Heegaard Splitting.}
        \label{f:ghs}
        \end{center}
        \end{figure}

We will sometimes depict a GHS schematically as in Figure \ref{f:ghs}. Often when we do this we will also need to represent compressing disks for some thick levels. Examples of this are the curved arcs depicted in the figure.

\begin{dfn}
\label{d:order}
For any surface $F$ let $c(F)=\sum \limits _n (2-\chi(F^n))^2$, where $\{F^n\}$ are the components of $F$. If $F_1$ and $F_2$ denote compact, embedded surfaces in a 3-manifold $M$ then we say $F_1 < F_2$ if $c(F_1) < c(F_2)$.
\end{dfn}

Note that this ordering is defined so that if $F_1$ is obtained from $F_2$ by a compression then $F_1 < F_2$. 

\begin{dfn}
Let $F^1=\{F^1_i\}$ and $F^2=\{F^2_j\}$ be two GHSs of a 3-manifold $M$. We say $F^1<F^2$ if $\{F^1_i\}_{i \ {\rm odd}}<\{F^2_j\}_{j \ {\rm odd}}$, where each set is put in non-increasing order and then the comparison is made lexicographically. 
\end{dfn}

\begin{dfn}
A GHS $\{F_i\}$ is said to be {\it strongly irreducible} if each thick level $F_i$ is strongly irreducible in the submanifold cobounded by $F_{i-1}$ and $F_{i+1}$. 
\end{dfn}

We now define two ways to get from a GHS which is not strongly irreducible to a smaller one. Suppose $S^*=\{F_i\}$ is a GHS. Suppose further that $D-E$ is an edge in $\Gamma (F_i)$, where $D$ and $E$ are disks in the submanifold of $M$ co-bounded by $F_{i-1}$ and $F_{i+1}$, for some odd $i$. Let $F_D$ denote the surface obtained from $F_i$ by compression along $D$, and $F_E$ denote the surface obtained from $F_i$ by compression along $E$. If $D \cap E =\emptyset$, then let $F_{DE}$ denote the surface obtained from $F_i$ by compression along both $D$ and $E$. There are now two cases, with several subcases:

\begin{enumerate}
    \item $D \cap E =\emptyset$
        \begin{enumerate}
            \item $F_D \ne F_{i-1}$, $F_E \ne F_{i+1}$. 
            \\ Remove $F_i$ from $S^*$. In it's place, insert $\{F_D, F_{DE}, F_E\}$ and reindex.  
            \item $F_D = F_{i-1}$, $F_E \ne F_{i+1}$. 
            \\ Replace $\{F_{i-1}, F_i\}$ with $\{F_{DE},F_E\}$ in $S^*$.
            \item $F_D \ne F_{i-1}$, $F_E = F_{i+1}$. 
            \\ Replace $\{F_{i}, F_{i+1}\}$ with $\{F_D,F_{DE}\}$.
            \item $F_D = F_{i-1}$, $F_E = F_{i+1}$. 
            \\ Replace $\{F_{i-1}, F_{i}, F_{i+1}\}$ with $F_{DE}$ and reindex. 
        \end{enumerate}
    \item $|D \cap E|=1$ (In this case $F_D$ and $F_E$ co-bound a product region of $M$)
        \begin{enumerate}
            \item $F_D \ne F_{i-1}$, $F_E \ne F_{i+1}$.
            \\ Replace $F_{i}$ in $S^*$ with $F_D$.
            \item $F_D = F_{i-1}$, $F_E \ne F_{i+1}$. 
            \\ Remove $\{F_{i-1}, F_{i}\}$ from $S^*$ and reindex. 
            \item $F_D \ne F_{i-1}$, $F_E = F_{i+1}$. 
            \\ Remove $\{F_{i}, F_{i+1}\}$ from $S^*$ and reindex. 
            \item $F_D = F_{i-1}$, $F_E = F_{i+1}$. 
            \\ Remove $\{F_{i}, F_{i+1}\}$ {\it or} $\{F_{i-1}, F_{i}\}$ from $S^*$ and reindex. 
        \end{enumerate}
\end{enumerate}

\begin{figure}[htbp]
        \psfrag{a}{(a)}
        \psfrag{b}{(b)}
        \vspace{0 in}
        \begin{center}
        \epsfxsize=4 in
        \epsfbox{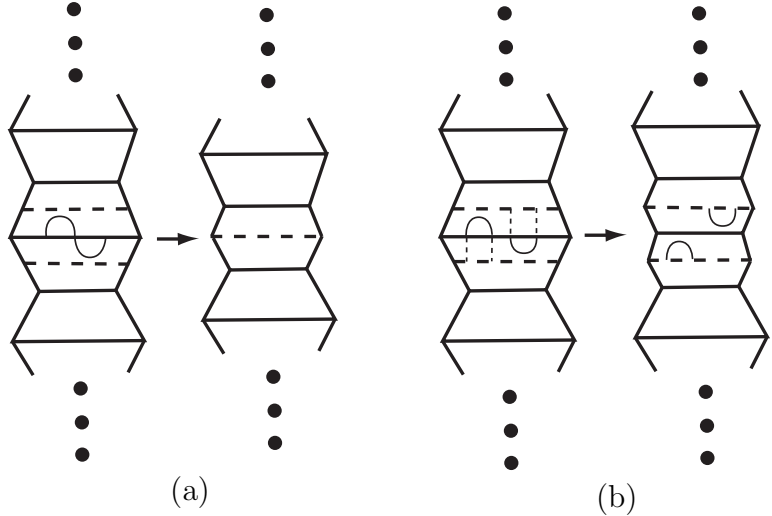}
        \caption{(a) A destabilization. (b) A weak reduction.}
        \label{f:reddfn}
        \end{center}
        \end{figure}

We leave it as an exercise to show that the new sequence thus defined is a GHS. In Case 1 ($D \cap E =\emptyset$) we say the new GHS was obtained from the old one by the {\it weak reduction} $D-E$. In Case 2 ($|D \cap E|=1$) we say the new GHS was obtained by a {\it destabilization}. Each of these operations is represented schematically in Figure \ref{f:reddfn}. Note that if the GHS $S^1$ is obtained from the GHS $S^2$ by weak reduction or destabilization then $S^1<S^2$. 

For readers unfamiliar with Generalized Heegaard Splittings we pause here for a moment to tie these concepts to more familiar ones. Suppose $M$ is a closed 3-manifold and $\{F_i\}_{i=0}^{2n}$ is a GHS of $M$. Since $M$ is closed we must have $F_0=F_{2n}=\emptyset$. If, in addition, $n=1$ then our GHS looks like $\{\emptyset, F_1, \emptyset\}$. By definition $F_1$ is a Heegaard splitting of $M$. If $\{\emptyset, F'_1, \emptyset\}$ was obtained from $\{\emptyset, F_1, \emptyset\}$ by a destabilization then the Heegaard splitting $F_1$ is a stabilization of the Heegaard splitting $F'_1$.

\begin{dfn}
A {\it Sequence Of GHSs} (SOG) of a 3-manifold is a sequence $\{F^j\}_{j=1}^n$ such that for each $k$ between 1 and $n-1$ one of the GHSs $F^k$ or $F^{k+1}$ is obtained from the other by a weak reduction or destabilization.
\end{dfn}

{\it Notation:} We will always use subscripts to denote surfaces, superscripts to denote GHSs, and a boldface font to denote an entire SOG. Hence, $F_i^j$ is the $i$th surface of $F^j$, which is the $j$th GHS of the SOG $\bf F$. If $\bf F$ is a SOG of a 3-manifold $M$ then $M^k_i$ will always denote the submanifold of $M$ cobounded by $F^k_{i-1}$ and $F^k_{i+1}$.

\begin{dfn}
If $\bf F$ is a SOG of $M$ and $k$ is such that $F^{k-1}$ and $F^{k+1}$ are both obtained from $F^k$ by weak reduction or destabilization then we say the GHS $F^k$ is {\it maximal} in $\bf F$. Similarly, if $k$ is such that $F^k$ is obtained from both $F^{k-1}$ and $F^{k+1}$ by weak reduction or destabilization then we say it is {\it minimal}.
\end{dfn}

\begin{dfn}
Let $M$ be a 3-manifold with two GHSs $F^-$ and $F^+$. We say $\{F^j\}_{j=1}^n$ is a SOG {\it from $F^-$ to $F^+$} if $F^1=F^-$ and $F^n=F^+$. 
\end{dfn}

In \cite{crit}, we prove the following Lemmas.

\begin{lem}
\label{l:oddstrirr}
If $F^-$ and $F^+$ are strongly irreducible GHSs of a 3-manifold $M$ then there is a SOG $\bf F$ from $F^-$ to $F^+$ such that if $F^k$ is maximal in $\bf F$ then there is exactly one thick level $F^k_i$ of $F^k$ which is not strongly irreducible in $M^k_i$.
\end{lem}

\begin{lem}
Among all SOGs described by the conclusion of Lemma \ref{l:oddstrirr} there is one such that if $F^k$ is maximal in $\bf F$ then there is exactly one thick level $F^k_p$ of $F^k$ which is critical in $M^k_p$.
\end{lem}

To prove this Lemma we establish the following stronger result which will be needed later:

\begin{lem}
\label{l:critmax}
Among all SOGs described by the conclusion of Lemma \ref{l:oddstrirr} there is one such that if $F^k$ is maximal in $\bf F$, $F^{k-1}$ is obtained from $F^k$ by the weak reduction or destabilization $D-E \in \Gamma(F^k_p)$, and $F^{k+1}$ is obtained from $F^k$ by the weak reduction or destabilization $D'-E'$ then $D-E$ is in a different component of $\Gamma (F^k_p)$ than $D'-E'$. 
\end{lem}

\begin{lem}
\label{l:critsfce}
If $\bf F$ is as described by the conclusions of Lemmas \ref{l:oddstrirr} and \ref{l:critmax} and $F^k$ is maximal in $\bf F$ then every thin level for $F^k$ is incompressible in $M$. In particular, if $M$ is non-Haken then $F^k=\{\emptyset, F^k_1, \emptyset\}$ for some critical Heegaard surface $F^k_1$ of $M$.
\end{lem}

\begin{lem}
\label{l:minHeegaard}
Among all SOGs described by the conclusion of Lemma \ref{l:critsfce} there is one such that if $F^k$ is minimal in $\bf F$ then each thick level $F^k_i$ is strongly irreducible in $M^k_i$ and each thin level of $F^k$ is incompressible in $M$. In particular, if $M$ is non-Haken then $F^k=\{\emptyset, F^k_1, \emptyset\}$ for some strongly irreducible Heegaard surface $F^k_1$ of $M$. 
\end{lem}

\begin{lem}
If $F^k=\{\emptyset, F^k_1,\emptyset\}$ is maximal in $\mathcal S$ and $F^l=\{\emptyset, F^l_1,\emptyset\}$ is the next (or previous) minimal GHS in $\mathcal S$ then $F^k_1$ is a stabilization of $F^l_1$.
\end{lem}

Combining the above Lemmas leads to: 

\begin{thm}
\label{t:critmain}
Suppose $F$ and $F'$ are non-trivial strongly irreducible Heegaard splittings of a non-Haken 3-manifold $M$. Then there is a SOG of $M$ from $\{\emptyset, F, \emptyset\}$ to $\{\emptyset, F', \emptyset\}$ with a unique maximal GHS $F^k=\{\emptyset, F^k_1, \emptyset\}$, where $F^k_1$ is a critical Heegaard splitting of $M$ and is the minimal genus common stabilization of $F$ and $F'$. 
\end{thm}

\subsection{GHSs of pairs $(M,\Sigma)$}

Now that we have seen how critical surfaces naturally arise we must strive to understand relatively critical surfaces. In this section we will review the relevant definitions and results from \cite{machine} pertaining to GHSs of pairs $(M,K)$, where $K$ is a properly embedded 1-manifold, and pairs $(M,\Gamma)$, where $\Gamma$ is a graph. In the next section we will look at SOGs of pairs. Only then will we be prepared to see how relatively critical surfaces naturally arise in various contexts.

To proceed we must first understand arcs in a compression body. If $W$ is a compression body recall that $W$ can be built by starting with a product $F \times I$ and attaching 2- and 3-handles to $F \times \{1\}$. Anything that remains of $F \times \{1\}$ after the attachment becomes part of $\partial _- W$. We say an arc $k$ is {\it straight} in $W$ if $k=\{p\} \times I$, where $p \in F$ is a point such that $\{p\} \times \{1\} \in \partial _-W$. 

We are now ready to generalize the definition of a compression body:

\begin{dfn}
A {\it $K$-compression body} $(W;K)$ is 
\begin{enumerate}
    \item A 3-manifold $W$ which can be obtained by starting with some surface $F$ (not necessarily connected), forming the product $F \times I$, attaching some number of 2-handles to $F \times \{1\}$, and capping some (but not necessarily all!) remaining 2-sphere boundary components with 3-balls. The boundary component $F \times \{0\}$ is $\partial _+ W$ and the other boundary component is $\partial _- W$.

    \item  A 1-manifold $(K, \partial K) \subset (W, \partial W)$ such that
    \begin{enumerate}
        \item $K$ is a disjoint union of embedded arcs
        \item each arc of $K$ has at least one endpoint on $\partial _+W$
        \item if $k$ is an arc of $K$ with $\partial k \subset \partial _+W$, then there is a disk, $D \subset W$, with $\partial D=k \cup \alpha$, $D \cap K=k$, and $D \cap \partial _+ W=\alpha$.  
        \item if $k$ is an arc of $K$ with one endpoint on $\partial _+W$, then $k$ is straight.
        \item For each 2-sphere component, $S$, of $\partial _-W$, $S \cap K \ne \emptyset$.
    \end{enumerate}
\end{enumerate}
A $K$-compression body $(W;K)$ is {\it non-trivial} if either $W$ is not a product or at least one arc of $K$ is not straight.
\end{dfn}

\begin{dfn}
If $K$ is a 1-manifold which is properly embedded in a 3-manifold $M$ then a {\it Heegaard splitting of the pair} $(M,K)$ is an expression of $M$ as a union of $K_i$-compression bodies $(W_1;K_1)$ and $(W_2;K_2)$ such that $\partial _+W_1=\partial _+ W_2$ and $K=K_1 \cup K_2$. Such a splitting is {\it non-trivial} if both $(W_1;K_1)$ and $(W_2;K_2)$ are non-trivial. 
\end{dfn}

When the context is clear, we will refer to the surface $\partial _+ W$ as a {\it Heegaard surface} of $(M,K)$. Given this definition, the following is now an immediate Corollary to Lemma \ref{l:index2Gabai}:

\begin{cor}
\label{c:incompboundary}
If $F$ is a critical Heegaard surface of $(M,K)$ then $(\partial M)^K$ is incompressible in $M^K$.
\end{cor}

We are now prepared to work our way through the definitions of GHS, weak reduction, destabilization, and SOG, and give corresponding definitions for each for pairs. 

\begin{dfn} 
Let $K$ be a 1-manifold which is properly embedded in a 3-manifold $M$. Let $\{G_j\}_{j=0}^{2n}$ be a sequence of closed, embedded, pairwise disjoint surfaces in $M$ which are transverse to $K$. Let $M_j$ denote the submanifold of $M$ co-bounded by $G_{j-1}$ and $G_{j+1}$. Then $\{G_j\}_{j=0}^{2n}$ is a {\it Generalized Heegaard Splitting (GHS) of the pair} $(M,K)$ if for each odd $j$ the surface $G_j$ is a non-trivial Heegaard splitting, or a union of non-trivial Heegaard splittings, of $(M_j, K)$ and $\partial M=G_0 \amalg G_{2n}$. 
\end{dfn}

\begin{dfn} 
If $\{G_j\}_{j=0}^{2n}$ is a GHS of $(M,K)$ then for each odd $j$ the surface $G_j$ is referred to as a {\it thick level} and for each even $j$ the surface $G_j$ is a {\it thin level}. 
\end{dfn}

\begin{ex}
\label{e:thin}
Suppose that $K \subset S^3$ is an arbitrary knot or link with no trivial components and $h$ is some standard height function on $S^3$ (so that for each $p \in (0,1)$, $h^{-1}(p)$ is a 2-sphere) which is a Morse function when restricted to $K$. Let $\{ q'_j \}$ denote the critical values of $h$ restricted to $K$ and let $q_j$ be some point in the interval $(q'_j, q'_{j+1})$. The following terminology is standard in {\it thin position} arguments (see \cite{gabai:87}).

If $j$ is such that $|K \cap h^{-1}(q_j)|>|K \cap h^{-1}(q_{j-1})|$ and $|K \cap h^{-1}(q_j)|>|K \cap h^{-1}(q_{j+1})|$ then we say the surface $h^{-1}(q_j)$ is a {\it thick level} of $K$. Similarly, if $|K \cap h^{-1}(q_j)|<|K \cap h^{-1}(q_{j-1})|$ and $|K \cap h^{-1}(q_j)|<|K \cap h^{-1}(q_{j+1})|$ then we say the surface $h^{-1}(q_j)$ is a {\it thin level} of $K$. 

Suppose there are $n$ thick levels for $K$. Let $G_0=G_{2n}=\emptyset$, $\{G_{2j-1}\}_{j=1}^n$ denote the set of thick levels of $K$, and $\{G_{2j}\}_{j=0}^{n-1}$ denote the set of thin levels. Then $\{G_j\}_{j=0}^{2n}$ is a GHS of  $(S^3,K)$. The thick and thin levels of this GHS are precisely the thick and thin levels of $K$.  
\end{ex}

We will sometimes depict a GHS of a pair schematically using the same type of figure as before (see Figure \ref{f:ghs}). However, now the curved arcs in the figure may represent relative compressing disks as well as compressing disks.

We now present one further generalization of the concept of a GHS. If $\Sigma$ is any properly embedded graph then let $\Sigma ^0$ and $\Sigma ^1$ denote the sets of vertices and edges of $\Sigma$. If we are discussing a particular graph, $\Sigma$, embedded in a 3-manifold, $M$, then $M^*$ will denote $M$ with a regular neighborhood of $\Sigma ^0$ removed. 

\begin{dfn}
Let $\Sigma$ be a properly embedded graph in a 3-manifold $M$. We define a {\it GHS of $(M,\Sigma)$} to be a GHS of $(M^*, \Sigma ^1)$.
\end{dfn}

Given a GHS of a pair $(M,\Sigma)$ we can define a GHS of $M$ by ``forgetting" $\Sigma$. 

\begin{dfn}
\label{d:underlyingGHS}
If $\{G_j\}_{j=0}^n$ is a GHS of $(M,\Sigma)$, where $M$ is an irreducible 3-manifold other than $B^3$ or $S^3$, then we define its {\it underlying GHS} $[\{G_j\}]$ as follows: 
\begin{enumerate}
    \item Let $s:\{0,...,n\} \rightarrow \{0,...,m\}$ be the onto, monotone function such that $s(i)=s(j)$ iff the submanifold of $M$ co-bounded by the non-$S^2$ components of $G_i$ and $G_j$ is a product. 
    \item For each $i$ between 0 and $m$ choose some $j \in s^{-1}(i)$ and let $F_i$ denote the non-$S^2$ components of $G_j$.  
    \item Let $\sigma$ be the maximal subset of $\{0,...,m\}$ such that $\{F_i\}_{i \in \sigma}$ is a GHS of $M$. 
    \item Define $[\{G_j\}]=\{F_i\}_{i \in \sigma}$.
\end{enumerate}
We leave it to the reader to check that $[\{G_j\}]$ is well defined up to isotopy.
\end{dfn}

The following Lemma is proved in \cite{machine}.

\begin{lem}
\label{l:underlyingodd}
Let $\Sigma$ be a properly embedded graph in an irreducible 3-manifold $M$ other than $S^3$ or $B^3$. Let $\{G_j\}$ be a GHS of $(M,\Sigma)$ and suppose $\{F_i\}=[\{G_j\}]$. Then for each thick (thin) level $F_p$ of $\{F_i\}$ there is a thick (thin) level of $\{G_j\}$ which becomes parallel to $F_p$ when all $S^2$ components are removed. 
\end{lem}

\begin{dfn}
Let $\Sigma$ be a properly embedded graph in a 3-manifold $M$. For any surface $F \subset M$ let $c(F;\Sigma)=\sum \limits _n \left(2-\chi(F^n_{\Sigma})\right)^2$, where $\{F^n\}$ are the components of $F$ and $F^n_{\Sigma}$ denotes $F^n$ with a neighborhood of $\Sigma$ removed. If $F_1$ and $F_2$ denote compact, embedded surfaces in $M$ 
then we say $F_1 <_{\Sigma} F_2$ if $c(F_1;\Sigma) < c(F_2;\Sigma)$.
\end{dfn}

\begin{dfn}
Let $\Sigma$ be a properly embedded graph in a 3-manifold $M$.  Let $G^1=\{G^1_i\}$ and $G^2=\{G^2_j\}$ be two GHSs of $(M,\Sigma)$. We say $G^1<_{\Sigma} G^2$ if $\{G^1_i|i\ {\rm odd}\}<_{\Sigma} \{G^2_j|j\ {\rm odd}\}$, where each set is put in non-increasing order and then the comparison is made lexicographically. 
\end{dfn}

Suppose $\{G_i\}$ is a GHS of $(M,\Sigma)$. Suppose further that for some odd $i$, $D-E$ is an edge in $\Sigma (G_i;K)$ such that both $D$ and $E$ lie in the submanifold of $M$ co-bounded by $G_{i-1}$ and $G_{i+1}$. Let $G_D$ denote the surface obtained from $G_{i}$ by compression along $D$ (if $D$ is a compressing disk) or by an isotopy guided by $D$ (if $D$ is a relative compression). Similarly, let $G_E$ denote the surface obtained from $G_{i}$ by compression along $E$ or by an isotopy guided by $E$. If $D \cap E =\emptyset$ then let $G_{DE}$ denote the surface obtained from $G_{i}$ by compression along (or isotopies guided by) both $D$ and $E$. The definitions of {\it weak reduction} and {\it destabilization} are now {\it exactly} the same as in section \ref{s:critreview}. 

\begin{ex}
\label{e:thin2}
Let $K$, $h$, and $\{q_j'\}$ be as in Example \ref{e:thin}. The {\it width} of $K$ is defined to be the quantity $\sum \limits _j |K \cap h^{-1}(q_j')|$. A knot is said to be in {\it thin position} if $h$ is chosen so that the width of $K$ is minimal (see \cite{gabai:87}).

Recall from Example \ref{e:thin} that our choice of $h$ induces a relative GHS $\{G_i\}$ of $(S^3;K)$. In thin position arguments a red relative compressing disk is often referred to as a ``strict high disk" and a blue relative compressing disk is a ``strict low disk". Most of the arguments that use thin position begin by assuming the opposite of what is to be proved about $K$ {\it and} that $K$ is in thin position. The contradiction is almost always the production of disjoint high and low disks or a high disk that meets a low disk in a point (in either case it is easy to show that the width of $K$ was not minimal). In our terminology, we would say that there was a weak reduction or a destabilization for $\{G_i\}$, and therefore there exists a smaller GHS.
\end{ex}

Recall that if $\Sigma$ is a properly embedded graph in a 3-manifold $M$ then $M^*$ denotes $M$ with a regular neighborhood of $\Sigma ^0$ removed.

\begin{dfn}
Let $\Sigma$ be a properly embedded graph in a 3-manifold $M$. Let $F$ be a closed, embedded, separating surface in $M^*$. Let $D$ be a $\Sigma ^1$-compression for $F$. If $D$ is a compressing disk for $F$ then we say it is {\it opaque}. Otherwise it is {\it transparent}.
\end{dfn}

Note that if $D$ is transparent disk then there are two cases. The first is that it is a relative compressing disk. The second is that it is a compressing disk for $F^{\Sigma ^1}$ which is not a compressing disk for $F$. In the latter case $\partial D$ bounds a disk on $F$ which meets $\Sigma^1$.

\begin{lem}
\label{l:sameGHS}
Suppose $G^2$ is a GHS of $(M,\Sigma)$ obtained from $G^1$ by the weak reduction or destabilization $D-E$. Then $[G^1] = [G^2]$ if and only if $D$ or $E$ is transparent.
\end{lem}

\begin{proof}
The proof is simply a matter of chasing through the definitions of weak reduction and destabilization. As this definition is given in eight different cases the full proof is quite lengthly and tedious. We will only present the first case, and leave it to the reader to check the rest. 

Recall from the definitions of weak reduction and destabilization that $D-E$ is an edge of $\Sigma (G_i;K)$, for some odd $i$. Recall also the surfaces $G_D$, $G_E$, and $G_{DE}$. We will prove Lemma \ref{l:sameGHS} in the case when $D \cap E=\emptyset$ and $G_D \ne G_{i-1}$, $G_E \ne G_{i+1}$. 

In this case $D-E$ is a weak reduction. To obtain $G^2$ from $G^1$ we were to remove $G_{i}$ from $G^1$ and in its place insert $\{G_D, G_{DE}, G_E\}$. Denote the non-$S^2$ components of a surface $G$ by $\bar G$. Now, if $D$ is a transparent disk then $\bar G_D$ is isotopic to  $\bar G_i$. Hence there is a compression body between $\bar G_D$ and $\bar G_{i+1}$. To form the underlying GHS of $G^2$ we include only those surfaces of $G^2$ that define a GHS. Hence, we skip $G_E$ and $G_{DE}$. That is, although $G^2$ contains the subsequence $\{G_{i-1},G_D, G_{DE}, G_E,G_{i+1}\}$, the underlying GHS of $G^2$ only contains the subsequence $\{\bar G_{i-1},\bar G_D,\bar G_{i+1}\}$, which is isotopic to $\{\bar G_{i-1},\bar G_{i},\bar G_{i+1}\}$. 
\end{proof}






\subsection{SOGs of pairs $(M,\Sigma)$} 

The next two definitions should be obvious to the reader who has followed things this far.

\begin{dfn}
Let $\Sigma$ be a graph properly embedded in a 3-manifold $M$. A {\it Sequence of GHSs (SOG) of $(M,\Sigma)$} is a sequence $\{G^j\}_{j=1}^n$ such that for each $k$ between 1 and $n-1$ one of the GHSs $G^k$ or $G^{k+1}$ is obtained from the other by a weak reduction or destabilization.
\end{dfn}

\begin{dfn}
If $\bf G$ is a SOG of $(M,K)$ and $k$ is such that $G^{k-1}$ and $G^{k+1}$ are both obtained from $G^k$ by weak reduction or destabilization then we say $G^k$ is {\it maximal} in $\bf G$. Similarly, if $k$ is such that $G^k$ is obtained from both $G^{k-1}$ and $G^{k+1}$ by weak reduction or destabilization then we say it is {\it minimal}.
\end{dfn}

Just as we defined the underlying GHS of a GHS of $(M,\Sigma)$ by ``forgetting" $\Sigma$, so to can we define an underlying SOG of a SOG of $(M,\Sigma)$.

\begin{dfn}
\label{d:underlyingSOG}
If ${\bf G}=\{G^j\}_{j=1}^n$ is a SOG of $(M,K)$, where $M$ is an irreducible 3-manifold other than $B^3$ or $S^3$, then we define its {\it underlying SOG} $[\bf G]$ as follows:
\begin{enumerate}
    \item Let $m$ be the smallest integer such that there is an onto, monotone function $s:\{1,...,n\} \rightarrow \{1,...,m\}$ such that $s(i)=s(j)$ iff $[G^i]=[G^j]$. 
    \item For each $i$ between 1 and $m$ choose some $j \in s^{-1}(i)$ and let $F^i=[G^j]$.
    \item Define $[{\bf G}]=\{F^i\}_{i=1}^m$.
\end{enumerate}
\end{dfn}


\begin{lem}
\label{l:underlyingGHS}
Suppose $\bf G$ is a SOG of $(M,\Sigma)$ and $\bf F=[\bf G]$. Then for each maximal (minimal) GHS $F^k$ of $\bf F$ there is a maximal (minimal) GHS of $\bf G$ whose underlying GHS is $F^k$.
\end{lem}

\begin{proof}
If $k$ is such that $F^k$ is maximal in $\bf F$ then $F^k> F^{k-1}$ and $F^k> F^{k+1}$. Let $a$ and $b$ be the smallest and largest elements of the set $\{j|[G^j]=F^k\}$. Then $[G^{a-1}]=F^{k-1}$, $[G^{b+1}]=F^{k+1}$ and $[G^c]=F^k$, for any number $c$ such that $a \le c \le b$. $F^k> F^{k-1}$ then implies that $G^{a}>_{\Sigma} G^{a-1}$. Let $c$ denote the largest integer between $a-1$ and $b+1$ (inclusive) such that $G^{c+1}>_{\Sigma} G^c$. Since $F^k> F^{k+1}$ implies $G^{b}>_{\Sigma} G^{b+1}$ we know $c<b$. But then we have $G^{c+1}>_{\Sigma} G^c$ and $G^{c+1}>_{\Sigma} G^{c+2}$. Hence, $G^{c+1}$ is maximal in $\bf G$.

The proof in the case that $F^k$ is minimal in $\bf F$ is simply a matter of switching all of the inequalities.
\end{proof}

\subsection{Comparing SOGs}

\begin{dfn}
Suppose $G$ is a relative Heegaard splitting of $(M,\Sigma)$. We say a path in $\Gamma (F;\Sigma)$ is {\it expanded} if its vertices have been realized by a sequence of disks $\{C_j\}$ such that for all $j$, 
\begin{enumerate}
    \item $C_j$ and $C_{j+1}$ are on opposite sides of $F$ and
    \item $C_j \cap C_{j+2}=\emptyset$.
\end{enumerate}
\end{dfn}

In the proof of Lemma \ref{l:critmax} we show that if $D-E$ and $D'-E'$ are edges in the same component of $\Gamma (F_i;K)$ then there is an expanded path which connects them. Although this proof is written for non-relative surfaces it works verbatim for relative ones as well.

We now define a function $\delta$ on maximal GHSs of $(M,\Sigma)$. 

\begin{dfn}
Suppose $G^k$ is a maximal GHS of a SOG $\bf G$ of $(M,\Sigma)$. Suppose further that $G^{k-1}$ is obtained from $G^k$ by the weak reduction or destabilization $D-E \in \Gamma(G^k_p;\Sigma^1)$ and $G^{k+1}$ is obtained from $G^k$ by the weak reduction or destabilization $D'-E' \in \Gamma(G^k_q;\Sigma^1)$. If $p=q$ and there is a path from $D-E$ to $D'-E'$ in $\Gamma(G^k_p;\Sigma^1)$ then define $\delta _{\bf G} (G^k)$ to be the length of the shortest expanded path. If $p=q$ and there is no such path then define $\delta _{\bf G} (G^k)=\infty$. Finally, if $p \ne q$ then define $\delta _{\bf G} (G^k)=1$.
\end{dfn}

{\bf Note.} By ``$\infty$" we simply mean some formal symbol for which the statement $\infty > n$ is true for all $n \in \mathbb Z$.

\begin{dfn}
If $H^*$ and $\bf G$ are a GHS and an SOG of $(M,\Sigma)$ then let $\sigma _{H^*}({\bf G})=\{j|G^j$ is maximal in ${\bf G}$ and $G^j=H^*\}$. The {\it multiplicity set of $H^*$ in $\bf G$}, $m_{H^*}(\bf G)$, is then defined to be the ordered set $\{\delta _{\bf G} (G^j)|j \in \sigma _{H^*}({\bf G})\}$, where we include repetitions and order in non-increasing order. 
\end{dfn} 

Note that if $G^k$ is a maximal GHS of a SOG $\bf G$ of $(M,\Sigma)$ then $|m_{G^k}(\bf G)|$ is precisely the number of times that a maximal GHS of $\bf G$ is equal to $G^k$. This is the justification for the choice of the term ``multiplicity set" for $m_{G^k}(\bf G)$.

\begin{dfn}
\label{d:orderSOG}
Let $\bf G$ and $\bf H$ be two SOGs of $(M,\Sigma)$. We say $\bf H<\bf G$ if there is a GHS $G^*$ which is maximal in $\bf G$ such that 
\begin{enumerate}
    \item $m_{G^*}({\bf H})$ is smaller than $m_{G^*}({\bf G})$ (were the comparison is made lexicographically), and
    \item for each GHS $E^*>G^*$ which is maximal in either $\bf G$ or $\bf H$, $m_{E^*}({\bf H})=m_{E^*}({\bf G})$.
\end{enumerate}
If neither $\bf G<\bf H$ nor $\bf H<\bf G$ is true then we say $\bf G \sim \bf H$.
\end{dfn}

In practice we will only be using this definition to show that some operation performed on a SOG yields one which is smaller. In all cases such an operation will begin with a maximal GHS $G^k$ of some SOG $\bf G$. The result will be a new SOG $\bf H$ in which $m_{G^k}({\bf H})<m_{G^k}({\bf G})$ (often $G^k$ will not appear at all as a maximal GHS of $\bf H$, in which case $m_{G^k}({\bf H})=\emptyset$). As no maximal GHS which is larger than $G^k$ will ever be affected in the transition from $\bf G$ to $\bf H$ Definition \ref{d:orderSOG} says that ${\bf H}<{\bf G}$.

\begin{dfn}
\label{d:dilute}
If $\Phi$ is a set of SOGs of $(M,\Sigma)$ and ${\bf G} \in \Phi$ is such that for each ${\bf H} \in \Phi$ either $\bf G< \bf H$ or $\bf G \sim \bf H$ then we say $\bf G$ is a {\it dilute} SOG in $\Phi$. When $\Phi$ is understood we simply say $\bf G$ is {\it dilute}.
\end{dfn}

\subsection{Dilute SOGs of $(B^3,T^1)$.}

The results of this section will not be used in the remainder of the paper. However, they are important on their own as they serve to illustrate the kind of situation in which 2-normal surfaces arise. The proofs also serve as a good warm-up for the section which follows.

Let $M$ be a 3-manifold with traingulation $T$ and $B$ an embedded 3-ball in $M$ with normal boundary. Suppose $B \cap T^0=\emptyset$ and the only normal 2-spheres in $B$ are copies of $\partial B$. Then Rubinstein \cite{rubinstein:93}, Thompson \cite{thompson:94}, and the author \cite{machine} have shown that $B$ contains an almost normal 2-sphere with an octagon. In this section we prove the following:

\begin{thm}
\label{t:almostnormalpair}
If $B$ contains two almost normal 2-spheres, each with an octagon, then $B$ contains a 2-normal 2-sphere with either two octagons or one 12-gon. 
\end{thm}

\begin{proof}
Let $S$ and $S'$ be almost normal 2-spheres in $B$ with octagons. Then $\{\partial B, S,\emptyset\}$ and $\{\partial B, S',\emptyset\}$ are strongly irreducible GHSs of $(B,T^1)$. Let $\Phi$ denote the set of all SOGs of $(B,T^1)$ from $\{\partial B, S,\emptyset\}$ to $\{\partial B, S',\emptyset\}$ such that for all ${\bf F} \in \Phi$ and all $i$ and $j$ every component of the surface $F_i^j$ is homeomorphic to $S^2$ or is $\emptyset$. Let $\bf G$ denote a dilute SOG in $\Phi$. 

Recall that for a SOG such as $\bf G$ the notation $M^j_i$ denotes the submanifold of $M$ which lies between $G^j_{i-1}$ and $G^j_{i+1}$.

\begin{clm}
\label{c:reloddstrirr}
If $G^k$ is maximal in $\bf G$ then there is exactly one thick level $G^k_i$ of $G^k$ which is not a strongly irreducible Heegaard splitting of $(M^k_i,T^1)$.
\end{clm}

\begin{proof}
Since $G^k$ is maximal in ${\bf G}$ there is some thick level $G^k_p$ such that $G^{k-1}$ is obtained from $G^k$ by a weak reduction or destabilization corresponding to an edge $D-E$ in $\Gamma (G^k_p;T^1)$. Similarly, there is a thick level $G^k_q$ such that $G^{k+1}$ is obtained from $G^k$ by a weak reduction or destabilization corresponding to an edge $D'-E'$ in $\Gamma (G^k_q;T^1)$. 

We first claim that $p=q$. If not then replace $G^k$ with $H^*$ in ${\bf G}$, where $H^*$ is the GHS of $(M^k_q,T^1)$ obtained from $G^{k-1}$ by the weak reduction or destabilization $D'-E'$. Since $H^*$ can also be obtained from $G^{k+1}$ by the weak reduction or destabilization $D-E$ our substitution has defined a new SOG ${\bf G}'$ of $(B,T^1)$. Note that $H^*$ is not maximal in ${\bf G}'$, so the multiplicity set of $G^k$ is smaller in ${\bf G}'$ than it is in ${\bf G}$.

We now show that for all odd $i \ne p$, $G^k_i$ is a strongly irreducible Heegaard splitting of $M^k_i$ relative to $T^1$. By way of contradiction, assume that $G^k_r$ is a thick level which is not strongly irreducible (where $r \ne p$) and let $D^*-E^*$ be an edge in $\Gamma (G^k_r;T^1)$. Let $H^-$, $H^0$, and $H^+$ denote the GHSs obtained from $G^{k-1}$, $G^k$, and $G^{k+1}$ by the weak reduction or destabilization $D^*-E^*$. Now replace $G^k$ in ${\bf G}$ with the subsequence $\{H^-,H^0,H^+\}$ to define a new SOG ${\bf G}''$ of $(B,T^1)$. As before, the multiplicity set of $G^k$ is smaller in ${\bf G}''$ than it is in ${\bf G}$. 
\end{proof}

\begin{clm}
If $\bf G$ is maximal in $\bf G$ then there is exactly one thick level $G^k_p$ of $G^k$ such that $G^k_p$ is critical in $M^k_p$ relative to $T^1$.
\end{clm}

\begin{proof}
Since $G^k$ is maximal in ${\bf G}$ there is some odd $p$ such that $G^{k-1}$ is obtained from $G^k$ by a weak reduction or destabilization corresponding to an edge $D-E$ in $\Gamma (G^k_p;T^1)$. Similarly, there is an odd $q$ such that $G^{k+1}$ is obtained from $G^k$ by a weak reduction or destabilization corresponding to an edge $D'-E'$ in $\Gamma (G^k_q;T^1)$. Claim \ref{c:reloddstrirr} implies $p=q$. It is our goal now to show that $D-E$ and $D'-E'$ lie in different components of $\Gamma (G^k_p;T^1)$. If this is not the case then note that $\delta _{\bf G}(G^k) \in \mathbb Z^+$. That is, there is some expanded path in $\Gamma (G^k_p;T^1)$ from $D-E$ to $D'-E'$ whose length is equal to $\delta _{\bf G} (G^k)$. 

If $\delta _{\bf G}(G^k)=1$ then $D-E$ is the same as $D'-E'$, so removal of the subsequence $\{G^k, G^{k+1}\}$ from ${\bf G}$ defines a new relative SOG ${\bf G}'$ of $(M,T^1)$ such that $m_{G^k}({\bf G}')<m_{G^k}({\bf G})$. 

If $\delta _{\bf G}(G^k) = 2$ then a shortest expanded path is of the form $D-E=E'-D'$ or $E-D=D'-E'$. Assume the former. There are now precisely three cases (up to symmetry) for the configuration of $D$, $E$, and $D'$, depending on whether $|D \cap E|=0$ or 1 and on whether $|D' \cap E|=0$ or 1. In each of these cases we can define a new SOG ${\bf G}'$ of $(B,T^1)$ where the multiplicity set of $G^k$ is smaller, exactly as in the proof of Lemma \ref{l:critmax}. 

Now suppose $\delta _{\bf G}(G^k) = n$ for some $n>2$. Choose some $m$ between 2 and $n-1$ and let $D^*-E^*$ denote the $m$th edge of the expanded path (so that $D^*-E^*$ is not equal to either $D-E$ to $D'-E'$). Let $G^*$ denote the GHS of $(M,T^1)$ obtained from $G^k$ by the weak reduction or destabilization corresponding to $D^*-E^*$. Now, let ${\bf G}'$ denote the SOG  obtained from ${\bf G}$ by inserting the subsequence $\{G^*, G^k\}$ just after $G^k$. Note that the maximal GHS $G^k$ appears one more time in ${\bf G}'$ than in ${\bf G}$. However, the set $m_{G^k}({\bf G}')$ can be obtained from the set $m_{G^k}({\bf G})$ by removing the number $n$ and inserting the numbers $m$ and $n-m+1$. Under the lexicographical ordering this is a decrease since both of these numbers are less than $n$. 

As we have ruled out all possibilities for $\delta _{\bf G}(G^k)$ other than $\infty$ we may now conclude that $D-E$ and $D'-E'$ are in different components of $\Gamma (G^k_p;T^1)$ and hence, $G^k_p$ is relatively critical. 
\end{proof}

Theorem \ref{t:classify} now implies that $G^k_p$ is a 2-normal surface in $B$. 

\begin{clm}
\label{c:notubes}
$G^k_p$ contains no tubes.
\end{clm}

\begin{proof}
Let $D$ be a compressing disk for a tube of $G^k_p$. Since $G^k_p$ is a union of 2-spheres $\partial D$ is separating on $G^k_p$. Hence, if $E$ is also a compressing disk then $D$ cannot meet $E$ in a point. We conclude that in all cases $D \cap E=\emptyset$, so that $D-E$ represents a weak reduction of $G^k$. Now, let $G^*$ denote the GHS obtained from $G^k$ by performing the weak reduction $D-E$. 

To perform the weak reduction $D-E$ we are to first form the surfaces $G_D$, $G_E$, and $G_{DE}$, obtained from $G^k_p$ by $T^1$-compression along $D$, $E$, and both $D$ and $E$. $G^*$ will then contain the subsequence $\{G_D, G_{DE}, G_E\}$, where $G_D$ and $G_E$ are thick levels. Note that $G_D$ contains a normal 2-sphere component. Hence, every GHS obtained from $G^*$ by a sequence of weak reductions and destabilizations will have a thick level with a normal 2-sphere component. 

Let $G'$ denote the strongly irreducible GHS of $(B,T^1)$ obtained from $G^*$ by performing as many weak reductions and destabilizations as possible. In \cite{machine} we show that the thin levels of $G'$ are normal and the thick levels are almost normal. As the only normal 2-spheres in $B$ are parallel to $\partial B$ we conclude that $G'=\{\partial B, Q, \emptyset\}$, where $Q$ is some almost normal surface. By the result of the preceding paragraph we know that $Q$ has some component that is a normal 2-sphere.

We now have an impossible situation. $Q$ must normalize to $\partial B$ on one side and $\emptyset$ on the other. But a normal component of $Q$ will not be affected by such a normalization. Hence after normalizing in either direction we can never obtain $\emptyset$.
\end{proof}

Claim \ref{c:notubes} implies that the exceptional pieces of $G^k_p$ are a pair of octagons or a single 12-gon. The only case in which Theorem \ref{t:almostnormalpair} can now be false is if $G^k_p$ contains two components, each with an octagon in a different tetrahedron. But then $B$ would contain disjoint, non-parallel almost normal 2-spheres $S_1$ and $S_2$. 

By \cite{jr:02} we know that almost normal surfaces act as ``barriers" to normalization. Hence normalizing $S_1$ to the side that contains $S_2$ must produce a non-empty normal surface which cannot be parallel to $\partial B$, a contradiction. We conclude then that $G^k_p$ contains a component with two octagons or a single 12-gon.
\end{proof}

The preceding proof serves as an excellent warm-up for the results of the next section.

\subsection{Dilute SOGs of $(M,\Sigma)$ with prescribed underlying SOGs}

\begin{dfn}
If $F^1$ and $F^2$ are GHSs such that $F^1$ can be obtained from $F^2$ by a (possibly empty) sequence of weak reductions and destabilizations then we write $F^2 \rhd F^1$. 
\end{dfn}

Note that $F^2 \rhd F^1$ implies $F^2 \ge F^1$. 

\begin{dfn}
Suppose $\bf F$ is a SOG. Define the {\it absolute maxima} of $\bf F$, $ABS(\bf F)$, to be $\{F^i \in \bf F|$ for all $j, F^j \rhd F^i$ implies $F^j=F^i\}$.
\end{dfn}

Let $F^-$ and $F^+$ denote two strongly irreducible GHSs of a 3-manifold $M$. Let $\bf F$ be a SOG of $M$ from $F^-$ to $F^+$ which satisfies the conclusions of Lemmas \ref{l:oddstrirr} and \ref{l:critmax}. Let $\Omega$ denote the set of all SOGs of $M$ from $F^-$ to $F^+$ such that for all ${\bf H} \in \Omega$, $ABS({\bf H}) \subset ABS({\bf F})$

As usual, let $\Sigma$ denote a properly embedded graph in $M$. Let $G^-$ and $G^+$ denote two strongly irreducible GHSs of $(M,\Sigma)$ such that $[G^-]=F^-$ and $[G^+]=F^+$. Finally, let $\Phi$ denote the set of all SOGs of $(M,\Sigma)$ from $G^-$ to $G^+$ whose underlying SOG is in $\Omega$. The remainder of this section concerns results pertaining to dilute SOGs in $\Phi$.

\begin{lem}
\label{l:reloddstrirr}
If $\bf G$ is a dilute SOG in $\Phi$ and $G^k$ is maximal in $\bf G$ then there is exactly one thick level $G^k_i$ of $G^k$ which is not a strongly irreducible in $M^k_i$ relative to $\Sigma^1$.
\end{lem}

\begin{proof}
Since $G^k$ is maximal in ${\bf G}$ there is some thick level $G^k_p$ such that $G^{k-1}$ is obtained from $G^k$ by a weak reduction or destabilization corresponding to an edge $D-E$ in $\Gamma (G^k_p,\Sigma)$. Similarly, there is a thick level $G^k_q$ such that $G^{k+1}$ is obtained from $G^k$ by a weak reduction or destabilization corresponding to an edge $D'-E'$ in $\Gamma (G^k_q,\Sigma)$. 

\begin{clm}
\label{c:p=q}
$p=q$
\end{clm}

\begin{proof}
Suppose not. As in the proof of Claim \ref{c:reloddstrirr}, replace $G^k$ with $H^*$ in ${\bf G}$, where $H^*$ is the GHS of $(M,\Sigma)$ obtained from $G^{k-1}$ by the weak reduction or destabilization $D'-E'$. Since $H^*$ can also be obtained from $G^{k+1}$ by the weak reduction or destabilization $D-E$ our substitution has defined a new SOG ${\bf G}'$ of $M$. Note that $H^*$ is not maximal in ${\bf G}'$, so the multiplicity set of $G^k$ is smaller in ${\bf G}'$ than it is in ${\bf G}$. In order for this to contradict our assumption that ${\bf G}$ was dilute in $\Phi$ we must also show that $[{\bf G}']\in \Omega$. The proof breaks down into two cases:

\noindent {\bf Case 1.} {\it $D$, $E$, $D'$ and $E'$ are all opaque.} Lemma \ref{l:sameGHS} then implies that $[G^{k-1}] \ne [G_k]$ and $[G^k] \ne [G^{k+1}]$. Hence, it must be the case that $[G^k]$ is maximal in $[{\bf G}]$. 

{\bf Subcase 1.1.} {\it $[G^k] \in ABS([{\bf G}])$}. Since $[{\bf G}] \in \Omega$, $[G^k] \in ABS({\bf F})$, and hence satisfies the conclusion of Lemma \ref{l:oddstrirr}. However, as $p \ne q$ this is a contradiction. 

{\bf Subcase 1.2.} {\it $[G^k] \notin ABS([{\bf G}])$}. Then there is some $j \ne k$ such that $[G^j]\rhd [G^k]$. But then $[G^j] \rhd [G^{k-1}]$ and $[G^j] \rhd [G^{k+1}]$. As these GHSs are the only potential absolute maxima of $[{\bf G}']$ that are not absolute maxima of $[{\bf G}]$ we have shown that $ABS([{\bf G}'])=ABS([{\bf G}])$. Hence, $[{\bf G}']\in \Omega$.

\noindent {\bf Case 2.} {\it At least one of $D$, $E$, $D'$ and $E'$ is transparent.} Assume, without loss of generality, that $D$ is transparent. Lemma \ref{l:sameGHS} then implies that $[G^k]=[G^{k-1}]$, so the underlying SOG of ${\bf G}$ contains the subsequence $\{..., [G^{k-1}], [G^{k+1}],...\}$. Lemma \ref{l:sameGHS} also implies $[H^*]=[G^{k+1}]$, so the underlying SOG of ${\bf G}'$ also contains the subsequence $\{..., [G^{k-1}], [G^{k+1}],...\}$. This leads us to conclude that $[{\bf G}]=[{\bf G}']$, implying $[{\bf G}'] \in \Omega$. 
\end{proof}

To complete the proof of Lemma \ref{l:reloddstrirr} we show that for all odd $i \ne p$, $G^k_i$ is a strongly irreducible Heegaard splitting of $(M^k_i,\Sigma)$. This is again similar to the proof of Claim \ref{c:reloddstrirr}. By way of contradiction, assume that $G^k_r$ is a thick level which is not strongly irreducible (where $r \ne p$) and let $D^*-E^*$ be an edge in $\Gamma (G^k_r;\Sigma)$. Let $H^-$, $H^0$, and $H^+$ denote the GHSs obtained from $G^{k-1}$, $G^k$, and $G^{k+1}$ by the weak reduction or destabilization $D^*-E^*$. Now replace $G^k$ in ${\bf G}$ with the subsequence $\{H^-,H^0,H^+\}$ to define a new SOG ${\bf G}''$ of $(M,\Sigma)$. As before, the multiplicity set of $G^k$ is smaller in ${\bf G}''$ than it is in ${\bf G}$. To claim that this is a contradiction, we must once again show that $[{\bf G}''] \in \Omega$. The proof of this breaks down into two cases:

\noindent {\bf Case 1.} {\it Either $D^*$ or $E^*$ is transparent.} Lemma \ref{l:sameGHS} then implies that $[H^-]=[G^{k-1}]$, $[H^0]=[G^k]$, and $[H^+]=[G^{k+1}]$. Hence, it is easily seen that $[{\bf G}'']=[{\bf G}]$, and we are done. 

\noindent {\bf Case 2.} {\it $D^*$ and $E^*$ are both opaque.}

{\bf Case 2.1.} {\it $D$ and $E$ are both opaque.} As $r \ne p$ and every absolute maxima of $[{\bf G}]$ satisfies Lemma \ref{l:oddstrirr} it must be the case that $[G^k]$ is not an absolute maxima of $[{\bf G}]$. Hence, there is some $j \ne k$ such that $[G^j] \rhd [G^k]$. But then $[G^j] \rhd [G^{k-1}]$, $[G^j] \rhd [G^{k+1}]$, and $[G^j] \rhd [H^0]$. As these GHSs are the only potential absolute maxima of $[\bf G'']$ that are not absolute maxima of $[\bf G]$ we have shown that $ABS([{\bf G}''])=ABS([{\bf G}])$, and hence $[{\bf G}'']\in \Omega$. 

{\bf Case 2.2.} {\it $D'$ and $E'$ are both opaque.} This case and the previous one are symmetric. 

{\bf Case 2.3.} {\it At least one of $D$ and $E$, and at least one of $D'$ and $E'$, is transparent.} Lemma \ref{l:sameGHS} then implies that $[H^-]=[H^0]=[H^+]$. Hence, the only potential absolute maxima of $[{\bf G}'']$ that are not absolute maxima of $[{\bf G}]$ are $[G^{k-1}]$ and $[G^{k+1}]$. Suppose $[G^{k-1}]$ is an absolute maximum of $[{\bf G}'']$. Lemma \ref{l:sameGHS} implies that $[G^{k-1}]=[G^k]$, so $[G^k]$ must be an absolute maximum of $[{\bf G}]$. $[{\bf G}]\in \Omega$ then implies $[G^k] \in ABS({\bf F})$. Hence, $[G^{k-1}] \in ABS({\bf F})$.
\end{proof}

\begin{lem}
\label{l:relcritmax}
If $\bf G$ is a dilute SOG in $\Phi$ and $G^k$ is maximal in $\bf G$ then there is exactly one thick level $G^k_p$ of $G^k$ such that $G^k_p$ is relatively critical in $M^k_p$.
\end{lem}

\begin{proof}
Since $G^k$ is maximal in ${\bf G}$ there is some odd $p$ such that $G^{k-1}$ is obtained from $G^k$ by a weak reduction or destabilization corresponding to an edge $D-E$ in $\Gamma (G^k_p;\Sigma)$. Similarly, there is an odd $q$ such that $G^{k+1}$ is obtained from $G^k$ by a weak reduction or destabilization corresponding to an edge $D'-E'$ in $\Gamma (G^k_q;\Sigma)$. Lemma \ref{l:reloddstrirr} implies $p=q$. It is our goal now to show that $D-E$ and $D'-E'$ lie in different components of $\Gamma (G^k_p;\Sigma)$. If this is not the case then, as in the proof of Theorem \ref{t:almostnormalpair}, $\delta _{\bf G}(G^k) \in \mathbb Z^+$. That is, there is some expanded path in $\Gamma (G^k_p;\Sigma)$ from $D-E$ to $D'-E'$ whose length is equal to $\delta _{\bf G} (G^k)$.

\begin{clm}
\label{c:dne1}
$\delta _{\bf G}(G^k) \ne 1$.
\end{clm}

\begin{proof}
If $\delta _{\bf G}(G^k)=1$ then $D-E$ is the same as $D'-E'$, so removal of the subsequence $\{G^k, G^{k+1}\}$ from ${\bf G}$ defines a new relative SOG ${\bf G}'$ of $(M,\Sigma)$ such that $m_{G^k}({\bf G}')<m_{G^k}({\bf G})$. To say that this contradicts minimality we must also establish that $[{\bf G}']\in \Omega$. There are two cases:

\noindent {\bf Case 1.} {\it $D$ and $E$ are opaque.} Then by Lemma \ref{l:sameGHS} $[G^k]$ is maximal in $[{\bf G}]$. There are now two subcases:

{\bf Case 1.1.} $[G^k] \in ABS([{\bf G}])$. Then, as $[{\bf G}] \in \Omega$, $[G^k] \in ABS({\bf F})$ and hence satisfies the conclusion of Lemma \ref{l:critmax}. This is an obvious contradiction.

{\bf Case 1.2.} $[G^k] \notin ABS([{\bf G}])$. Then there is a $j \ne k$ such that $G^j \rhd G^k$. The only potential maximal GHS of $[{\bf G}']$ which is not a maximal GHS of $[{\bf G}]$ is $G^{k-1}$. But $G^j \rhd G^k$ implies $G^j \rhd G^{k-1}$, and so $ABS([{\bf G}'])=ABS([{\bf G}])$.

\noindent {\bf Case 2.} {\it $D$ or $E$ is transparent.} Then Lemma \ref{l:sameGHS} implies $[G^{k-1}]=[G^k]$, and so $ABS([{\bf G}'])=ABS([{\bf G}])$.
\end{proof}

\begin{clm}
$\delta _{\bf G}(G^k) \ne 2$.
\end{clm}

\begin{proof}
Suppose not. Then a shortest expanded path is of the form $D-E=E'-D'$ or $E-D=D'-E'$. Assume the former. There are now precisely three cases (up to symmetry) for the configuration of $D$, $E$, and $D'$, depending on whether $|D \cap E|=0$ or 1 and on whether $|D' \cap E|=0$ or 1. In each of these cases we can define a new SOG ${\bf G}'$ of $(M,\Sigma)$ where the multiplicity set of $G^k$ is smaller, exactly as in the proof of Lemma \ref{l:critmax}. For a contradiction, we must again show that $[{\bf G}']\in \Omega$. In our construction of ${\bf G}'$ the only new maximal GHSs which we potentially introduce are $G^{k-1}$ and $G^k$. We now list the cases (up to symmetry):

\noindent {\bf Case 1.} {\it $D$, $E$ and $D'$ are opaque.} See Case 1 of Claim \ref{c:dne1}. 

\noindent {\bf Case 2.} {\it Either both $D$ and $D'$, or $E$, is transparent.} See Case 2 of Claim \ref{c:dne1}.

\noindent {\bf Case 3.} {\it $D$ and $E$ are opaque and $D'$ is transparent.} If $[G^{k+1}]$ is an absolute maximum of $[{\bf G}']$ then $[G^k]$ is an absolute maximum of $[{\bf G}]$, and so $[G^k] \in ABS({\bf F})$. Lemma \ref{l:sameGHS} implies that $[G^{k+1}]=[G^k]$, so in this case $[G^{k+1}] \in ABS({\bf F})$. Note that $[G^{k-1}]$ cannot be an absolute maximum of $[{\bf G}']$ since $[G^k] \rhd [G^{k-1}]$ and $[G^{k+1}]=[G^k]$ implies $[G^{k+1}] \rhd [G^{k-1}]$. 
\end{proof}

\begin{clm}
$\delta _{\bf G}(G^k) \ne n$ for $n>2$.
\end{clm}

\begin{proof}
Suppose $\delta _{\bf G}(G^k) = n$ for some $n>2$. Choose some $m$ between 2 and $n-1$ and let $D^*-E^*$ denote the $m$th edge of the expanded path (so that $D^*-E^*$ is not equal to either $D-E$ to $D'-E'$). Let $G^*$ denote the GHS of $(M,\Sigma)$ obtained from $G^k$ by the weak reduction or destabilization corresponding to $D^*-E^*$. Now, let ${\bf G}'$ denote the SOG  obtained from ${\bf G}$ by inserting the subsequence $\{G^*, G^k\}$ just after $G^k$. Note that the maximal GHS $G^k$ appears one more time in ${\bf G}'$ than in ${\bf G}$. However, the set $m_{G^k}({\bf G}')$ can be obtained from the set $m_{G^k}({\bf G})$ by removing the number $n$ and inserting the numbers $m$ and $n-m+1$. Under the lexicographical ordering this is a decrease since both of these numbers are less than $n$. In addition we have not changed the set of absolute maxima so the underlying SOG of ${\bf G}'$ is also in $\Omega$. 
\end{proof}

As we have ruled out all possibilities for $\delta _{\bf G}(G^k)$ other than $\infty$ we may now conclude that $D-E$ and $D'-E'$ are in different components of $\Gamma (G^k_p;\Sigma)$ and hence, $G^k_p$ is relatively critical. 
\end{proof}

The next two lemmas will not be used in this paper, but do give important insight into the structure of a dilute SOG.

\begin{lem}
\label{l:relcritsfce}
If $\bf G$ is a dilute SOG in $\Phi$ and $G^k$ is maximal in $\bf G$ then every thin level of $G^k$ is incompressible and relatively incompressible in $M^K$.
\end{lem}

\begin{proof}
The proof is almost identical to that of Lemma \ref{l:critsfce}. Suppose $c$ is a loop that bounds a compressing disk or an arc which cobounds a relative compressing disk for some surface $G^k_q$, where $q$ is even. If $c$ is a loop then let $C$ denote a disk in $M^K$ such that $\partial C=c$. If $c$ is an arc then let $C$ be a disk in $M$ such that $\partial C=\gamma \cup c$, where $K \cap C=\gamma$. In either case choose $C$ so that $|C \cap (\bigcup \limits _{{\rm even} \ i} G^k_i)|$ is minimal. 

It follows from  Lemma 6.6. of \cite{tnbglex}, Corollary \ref{c:incompboundary}, and Lemmas \ref{l:reloddstrirr} and \ref{l:relcritmax} that for each odd $i$ $\partial M^k_i$ is incompressible and relatively incompressible in $M^k_i$. Hence, $C$ cannot lie entirely in $M^k_{q-1}$ or $M^k_{q+1}$. We conclude then that there is some loop or arc of intersection of the interior of $C$ with $\bigcup \limits _{{\rm even} \ i} G^k_i$. Let $\alpha$ denote an innermost such loop. Let $C'$ denote the subdisk of $C$ bounded by $\alpha$. $C'$ lies in $M^k_p$ for some odd number $p$. As $\partial M^k_p$ is incompressible in $M^k_p$, $\alpha$ must bound a disk $A$ on $\partial M^k_p$. 

Now, let $\beta$ be an innermost loop of $C \cap A$ and let $A'$ be the subdisk of $A$ bounded by $\beta$. Then we can use $A'$ to surger $C$ and thereby obtain a new disk with the same boundary as $C$, contradicting our minimality assumption.

We conclude then that $c$ cannot be a loop, and if $c$ is an arc then $C$ contains no loops of intersection with $\bigcup \limits _{{\rm even} \ i} G^k_i$. Let $\delta$ then denote an arc of intersection which is outermost on $C$. $\delta$ and a subarc of $\alpha$ cobound a subdisk $C''$ of $C$. $C''$ is then a relative compressing disk for $\partial M^k_i$ for some $i$, a contradiction.
\end{proof}

\begin{lem}
\label{l:relminHeegaard}
There is a dilute SOG in $\Phi$ whose minimal GHSs are strongly irreducible relative to $\Sigma^1$. 
\end{lem}

\begin{proof}
First, let $\bf G$ denote any dilute SOG of $\Phi$. Suppose $G^k$ is a minimal GHS of ${\bf G}$ which is not strongly irreducible. Then there is some weak reduction or destabilization for $G^k$. Let $G^*$ denote the result of such an operation. We can now define a new SOG by inserting the subsequence $\{G^*,G^k\}$ just after $G^k$. The new SOG will have one minimal GHS which is smaller. The rest of the minimal and maximal GHSs will remain unchanged. We may now repeat this operation to eventually arrive at the desired SOG. 
\end{proof}

\subsection{Pulling the rabbit out of the hat}

In this section we combine all of our results to prove that minimal genus common stabilizations of pairs of strongly irreducible Heegaard splittings can be isotoped to some normal form. 

\begin{dfn}
If $M$ is an irreducible triangulated 3-manifold then $S \subset M$ is a {\it maximal normal 2-sphere} in $M$ if there is no other normal 2-sphere which bounds a ball in $M$ which contains $S$. 
\end{dfn}

If $M$ is not homeomorphic to $B^3$ or $S^3$ and has a unique maximal normal 2-sphere $S$ then we recall two basic facts about $S$ from \cite{rubinstein:93} or \cite{thompson:94}. First, $S$ bounds a ball $B$ in $M$ which contains all of the vertices of $T^0 \cap M$. Second, there are no almost normal 2-spheres in $M$ that are disjoint from $B$. The same proof shows that there can also be no 2-normal 2-sphere in $M$ that are disjoint from $B$.

\begin{thm}
\label{t:main2}
The minimal genus common stabilization of any pair of strongly irreducible Heegaard splittings of a non-Haken 3-manifold is isotopic to an almost normal or 2-normal surface in any triangulation. 
\end{thm}

\begin{proof}
Let $M$ be a non-Haken 3-manifold and let $T$ be a triangulation of $M$. As in \cite{rubinstein:93} let $S$ be a maximal normal 2-sphere. $S$ bounds a ball $B$ which contains all of the vertices of $T^0 \cap M$. Collapsing $B$ to a point $v$ then turns $T^1$ into a properly embedded graph $\Sigma$ whose unique vertex is $v$. Note that any surface which is disjoint from $v$ after the collapse can be identified with a surface which is disjoint from $B$ before the collapse.

Let $F$ and $F'$ be a pair of strongly irreducible Heegaard splittings of $M$. By Theorem \ref{t:critmain} there is a SOG $\bf F$ of $M$ from $\{\emptyset, F, \emptyset\}$ to $\{\emptyset, F', \emptyset\}$ with a unique maximal GHS $F^k=\{\emptyset, F^k_1, \emptyset\}$, where $F^k_1$ is a critical Heegaard splitting of $M$ and is the minimal genus common stabilization of $F$ and $F'$. In particular, $\bf F$ satisfies the conclusions of Lemmas \ref{l:oddstrirr} and \ref{l:critmax}. 

Choose any strongly irreducible GHSs $G^-$ and $G^+$ of $(M,\Sigma)$ such that $[G^-]=\{\emptyset, F, \emptyset\}$ and $[G^+]=\{\emptyset, F', \emptyset\}$ and define the set $\Phi$ as in the previous section. 

Let $\bf G$ denote a dilute SOG of $(M,\Sigma)$ in $\Phi$. By Lemma \ref{l:underlyingGHS} there is a maximal GHS $G^p$ of $\bf G$ such that $[G^p]=F^k$. By Lemma \ref{l:underlyingodd} there is an odd number $q$ such that the non-$S^2$ component of $G^p_q$ is parallel to $F^k_1$. By Lemmas \ref{l:reloddstrirr} and \ref{l:relcritmax} $G^p_q$ is either strongly irreducible or critical relative to $\Sigma^1$. But then $G^p_q$ is either strongly irreducible or critical relative to $T^1$. If $G^p_q$ is relatively strongly irreducible then it is isotopic to an almost normal surface by Theorem 8.8 of \cite{machine}. If $G^p_q$ is relatively critical then it is isotopic to a 2-normal surface by Theorem \ref{t:classify}. 

Finally, notice that if any component of $G^p_q \backslash T^2$ that is not a triangle or quadrilateral lies in a 2-sphere component then we have an almost normal or 2-normal 2-sphere. This is a contradiction, as there are no such surfaces disjoint from $B$. We conclude that the non-$S^2$ component of $G^p_q$ ({\it i.e.} the component parallel to $F^k_1$) must be almost normal or 2-normal.
\end{proof}

\bibliographystyle{alpha}

\end{document}